\documentclass{elsart}
\usepackage{ifpdf}
\usepackage{graphicx,natbib,amssymb,lineno}
\usepackage {amsmath, latexsym, amscd, amssymb}
\usepackage {graphicx}
\usepackage[latin1]{inputenc}
\usepackage [french]{babel}
\ifpdf
\usepackage[%
  pdftitle={Groupes d'isom\'etries permutant doublement transitivement un ensemble  de droites vectorielles},%
  pdfauthor={Lucas Vienne},%
  pdfsubject={The preprint document class elsart},%
  pdfkeywords={automorphisme,  groupe, graphe, graphe fortement r\'egulier, groupe doublement transitif, graphe de Paley},%
  pdfstartview=FitH,%
  bookmarks=true,%
  bookmarksopen=true,%
  breaklinks=true,%
  colorlinks=true,%
  linkcolor=blue,anchorcolor=blue,%
  citecolor=blue,filecolor=blue,%
  menucolor=blue,pagecolor=blue,%
  urlcolor=blue]{hyperref}
\else
\usepackage[%
  breaklinks=true,%
  colorlinks=true,%
  linkcolor=blue,anchorcolor=blue,%
  citecolor=blue,filecolor=blue,%
  menucolor=blue,pagecolor=blue,%
  urlcolor=blue]{hyperref}
\fi

\makeatletter
\def\elsartstyle{%
    \def\normalsize{\@setfontsize\normalsize\@xiipt{14.5}}
    \def\small{\@setfontsize\small\@xipt{13.6}}
    \let\footnotesize=\small
    \def\large{\@setfontsize\large\@xivpt{18}}
    \def\Large{\@setfontsize\Large\@xviipt{22}}
    \skip\@mpfootins = 18\p@ \@plus 2\p@
    \normalsize
}
\@ifundefined{square}{}{\let\Box\square}
\makeatother

\newcommand{\PPP}{\mathbb{P}}

\newcommand{\EE}{\mathcal{E}}

\newcommand{\GG}{\mathcal{G}}

\newcommand{\SSS}{\mathcal{S}}

\newcommand{\ov}{\overline}

\newcommand{\lraa}{\longleftrightarrow}

\newcommand{\Lra}{\Leftrightarrow}

\newcommand{\ssi}{\Longleftrightarrow}
\newcommand{\x}{\times}

\newcommand{\al}{\langle}
\newcommand{\ar}{\rangle}

\newcommand{\Aut}{\mathrm{Aut}}

\newcommand{\stab}{\mathrm{stab}}

\newcommand{\nlb}{\nolinebreak}

\newcommand{\npb}{\nopagebreak}

\newcommand{\mm}{m\^eme }
\newcommand{\elt}{\'el\'ement }
\newcommand{\elts }{\'el\'ements }
\newcommand{\si }{si et seulement si }
\newcommand{\cad }{c'est-\`a-dire }

\newcommand{\Dem }{{\it D\'emonstration }}

\newcommand{\ttt}{th\'eor\`eme }

\newcommand{\oo}{\circ }
\newcommand{\cc }{\c{c}}
\newcommand{\ii }{\"{\i}}

\newcommand{\ev }{espace vectoriel }

\newcommand{\ga }{\alpha }
\newcommand{\gb }{\beta }
\newcommand{\gc }{\gamma }
\newcommand{\gd }{\delta }

\newcommand{\gve }{\varepsilon }
\newcommand{\gvf }{\varphi }

\newcommand{\gl}{\lambda}
\newcommand{\gm}{\mu}
\newcommand{\go }{\omega }
\newcommand{\gs}{\sigma}

\newcommand{\gC }{\Gamma }

\newcommand{\gL}{\Lambda}

\newtheorem{deff}{D\'efinition}
\newtheorem{lemf}{Lemme}
\newtheorem{thmf}{Th\'eor\`eme }

\newtheorem{propf}{Proposition}

\begin{document}

\begin{frontmatter}
\title{Groupes d'isom\'etries permutant doublement transitivement un ensemble  de droites vectorielles}

\author{Lucas Vienne}
\address{Departement de math\'ematiques. Universit\'e d'Angers. France}

\ead{lucas.vienne@univ-angers.fr}

\begin{abstract}
Soit  $n\geq 3$ un entier, et $G$  un groupe fini d'isom\'etries d'un
espace euclidien $E$ de dimension finie agissant deux fois
transitivement sur un ensemble de droites vectorielles
$\GG=\{U_1,\ldots,U_n\}$. Alors $\GG$  est une gerbe \'equiangulaire
de droites,  ce qui signifie que, pour $1\leq i,j\leq n$, il existe
des g\'en\'erateurs  $u_i$ des droites $U_i$, une constante $c$ et des  coefficients $\gve_{i,j}$ dans $\{-1,+1\}$ tels que \\[0,1cm]
\vspace{0,2cm}\centerline{$\forall i,j, \  1\leq i,j \leq n,$ \quad
  $\|u_i\|=1$,\  et si \ $i\neq j$ \ alors \ $(u_i|u_j)=\gve_{i,j}.c$}
On associe \`a ce syst\`eme g\'en\'erateur $(u_1,\ldots,u_n)$ le graphe simple $\gC$ sur l'ensemble $X=\{1,\ldots,n\}$ pour lequel deux points distincts $i$ et $j$ sont li\'es lorsque $\gve_{i,j}=-1$. \\
Dans cet article on \'etudie les propri\'et\'es du graphe $(\gC,X)$ impliqu\'ees par la double transitivit\'e de $G$, puis on recherche des graphes remplissant ces conditions et les groupes d'isom\'etries qui leurs sont associ\'es. C'est notamment le cas des graphes de Paley.
 \end{abstract}

\end{frontmatter}
\section{Introduction} 
 Soit $n\geq 3$ un entier, $E$ un \ev r\'eel de dimension finie et $G$ un sous-groupe fini du groupe lin\'eaire $GL(E)$ qui permute deux fois transitivement un ensemble de droites vectorielles $\GG=\{U_1,\ldots,U_n\}$. On sait qu'il existe un produit scalaire $\gvf=( \ |\ )$ sur $E$ pour lequel $G$ appara\^it comme un sous-groupe d'isom\'etries du groupe orthogonal $O_\gvf(E)$, et la double transitivit\'e de $G$ sur $\GG$ nous montre que les droites $U_i$ ($1\leq i\leq n$) font deux \`a deux un angle constant ; nous dirons que $\GG$ est une {\it gerbe \'equiangulaire}. Choisissant des g\'en\'erateurs $u_i$ de norme 1 pour chaque droite $U_i$, il existe une constante $c$ et des coefficients $\gve_{i,j}$ tous pris dans $\{-1,+1\}$ ($1\leq i,j\leq n$) tels  que \\[0,1cm]
 \vspace{0,2cm}\centerline{$\forall i,j,\  1\leq i,j\leq n , \  (u_i|u_i)=1 $
   \ et \  si \ $i\neq j$,   \  $(u_i,|u_j)=\gve_{i,j}.c.$ } 
Chacun des $2^n$ choix de ces syst\`emes de g\'en\'erateurs $(u_1,\ldots,u_n)$ est
associ\'e \`a un graphe simple $\gC$ sur l'ensemble
$X=\{1,\ldots,n\}$, deux sommets distincts $i$ et $j$ de $X$  \'etant
li\'es \si $\gve_{i,j}=-1$.\\
Dans la premi\`ere partie de cet article nous d\'ecrivons  les
propri\'et\'es g\'eom\'etriques des graphes $\gC$ impliqu\'ees par la
double transitivit\'e du groupe $G$ agissant sur la gerbe $\GG$.  Ces graphes seront dits {\it extensibles}.\\
Dans la deuxi\`eme partie on construit des graphes extensibles. On
montre notamment que ce sont toujours des graphes fortement
r\'eguliers et que  les graphes de Paley $P(q)$ sont tous
extensibles.\\
 La troisi\`eme partie s'appuie sur les r\'esultats  d'un article
 pr\'ec\'edent  ([5]) pour montrer que toute repr\'esentation d'un graphe extensible
 sur une gerbe  \'equi\-angulaire est la somme d'une repr\'esentation
 nulle et d'une repr\'esen\-tation dite r\'eduite, ces derni\`eres
 \'etant naturellement associ\'ees aux valeurs propres de la matrice
 du graphe $\gC$. On applique ce r\'esultat aux graphes construits dans
 la deuxi\`eme partie pour d\'eterminer les gerbes \'equiangulaires et les groupes d'isom\'etries qui leurs sont associ\'es.

\section{R\'esultats pr\'eliminaires}

Un entier $n\geq 3$ \'etant choisi, la
matrice $\EE=(\gve_{i,j})$ d'un graphe simple
$\gC$ sur l'ensemble $X=\{1,\ldots,n\}$ est donn\'ee par son $(i,j)$-\`eme coefficient  $\gve_{i,j}$ qui
vaut\ $-1$ \  si $i$ et $j$ sont li\'es (not\'e  $i\sim j$) et $1$
dans le cas contraire.\\
Donnons nous  un espace quadratique $(E,\ga)$ \cad un espace
vectoriel  r\'eel $E$ de dimension finie muni d'une forme bilin\'eaire
sym\'etrique $\ga $. Une repr\'esentation de $(\gC,X)$ dans 
$(E,\ga)$ est une application $u$ de $X$ dans
 $E$ not\'ee $u: i\to u_i$ telle que pour deux constantes $(\go,
 c)$ on ait \\[0,1cm]
$(1)$
\vspace{0,2cm}\centerline{$\forall i,j \in X,$ \  $\ga(u_i,u_i)=\go $
   \ et \  si \ $i\neq j$,   \  $\ga(u_i,u_j)=\gve_{i,j}.c.$ } 
On dit que $\go$ et $c$ sont les {\it param\`etres} de $u$ et que la
matrice not\'ee $S_\gC(u)$ ou plus simplement $S(u)$ de coefficient g\'en\'eral $\ga(u_i,u_j)$ $(i,j\in X)$ est
la {\it matrice de la repr\'esentation} $u$. Comme elle ne d\'epend
que de $\gC$ et des param\`etres $\go$ et $c$, on la note aussi
$S(\go,c)$ ; remarquons que $\EE=S(1,1)$. On appelle respectivement
{\it rang} et  {\it degr\'e} de $u$  le rang  de $S(u)$ et  la dimension de $E$. Enfin l'ensemble $\GG(u)$ des droites vectorielles engendr\'ees par les vecteurs $u_i\ (i\in X)$ s'appelle la gerbe {\it \'equiangulaire} ou {\it isom\'etrique} associ\'ee \`a $u$. On dira aussi que  le graphe $\gC$ est repr\'esent\'e sur la gerbe $\GG(u)$.

\subsection{Matrices et graphes associ\'es} 
Sauf indication contraire, le $(i,j)$-\`eme coefficient d'une matrice $M$ est not\'e $M_{i,j}$.
Pour tout $\gs$ dans le groupe $\SSS_X$ des permutations de $X$,
notons $P_\gs$ sa matrice, dont le $(i,j)$-\`eme coefficient vaut $P_{\gs,i,j}=\gd_{i,\gs(j)}$
(o\`u $\gd$ est le symbole de Kronecker usuel), et pour toute matrice
$M$ de type $n\x n$ posons   \mbox{${^\gs M}=P_\gs.M.P_\gs^{-1}$,} de sorte
que pour deux indices $i$, $j$ arbitraires dans $X$ on a \vspace{0,1cm}\mbox{${^\gs
  M}_{i,j}=M_{\gs^{-1}(i),\gs^{-1}(j)}$.} En remarquant que pour deux
permutations $\gs$ et $\gd$ de $X$ on a \ ${^{\gs\gd}M}={^\gs(}{^\gd
  M})$, on voit que l'ensemble des permutations $\gs$ de
$X$ telles que ${^\gs M}=M$ est un sous-groupe de $\SSS_X$ dit {\it
  stabilisateur} de $M$ et not\'e $\stab(M)$. De plus
$^\gs \EE$ est la matrice  du graphe ${^\gs\gC}$, image de $\gC$ par
la permutation $\gs$ qui est donc dans le groupe $\Aut(\gC)$ des
automorphismes du graphe $\gC$  \si  $^\gs\EE=\EE$. Autrement dit,  $\stab(\EE)=\Aut(\gC)$.\\
Soient $(\gC,X)$ un graphe repr\'esent\'e par une gerbe
isom\'etrique \mbox{$\GG=\{U_i \  |\  i\in X\}$} dans un espace quadratique
$(E,\ga)$, et  des g\'en\'erateurs $u_i$  des droites $U_i$ (pour
$i\in X$) satisfaisant \`a la condition $(1)$. Choisissant une suite
$(\nu_1,\ldots,\nu_n)$ de coefficients dans $\{-1,1\}$, les
g\'en\'erateurs \mbox{$u'_i=\nu_iu_i$} des droites $U_i$ ($i\in X$) satisfont
encore \`a la condition $(1)$ si l'on  remplace la matrice \mbox{$\EE=(\gve_{i,j})$} par la matrice $\EE'=(\gve'_{i,j})$  dont le $(i,j)$-\`eme coefficient vaut $\gve'_{i,j}=\nu_i\nu_j\gve_{i,j}$. Ce changement de g\'en\'erateurs induit aussi une modification du graphe $\gC$ associ\'e \`a la gerbe $\GG$ qui doit \^etre remplac\'e par le graphe $\gC'$ de matrice $\EE'$. Cette reflexion  nous conduit \`a la 
\begin{deff}$\ $\\ 
Soient $\gC$ et $\gC'$ deux graphes sur le \mm ensemble $X$ de
sommets, de matrices respectives $\EE=(\gve_{i,j})$ et
$\EE'=(\gve'_{i,j})$. On dit que les matrices $\EE$ et $\EE'$, ou les
graphes $\gC$ et $\gC'$, sont $\mathrm{ associ\acute{e}s}$ s'il existe une suite
$(\nu_1,\ldots,\nu_n)$ de coefficients, tous pris dans $\{-1,1\}$ tels que \\[0.2cm]
$(2)$ \vspace{0,2cm}\centerline{$\forall (i, j)\in X\x X , \qquad \gve'_{i,j}=\nu_i\nu_j.\gve_{i,j}$}
\end{deff}
La proposition suivante  nous donne quelques informations sur cette relation.
\begin{propf}$\ $\\ 
$1$. La relation d'association entre graphes ou matrices est une relation d'\'equivalence.\\
$2$.  L'ensemble $G(\EE)=G(\gC)$ des permutations $\gs$ de $X$ telles
que les matrices $\EE$ et $^\gs\EE$, ou les graphes $\gC$ et $^\gs\gC$, soient associ\'es est  un sous-groupe de $\SSS_X$ contenant le groupe $\stab(\EE)=\Aut(\gC)$.\\
 $3$. Si deux graphes $(\gC,X)$ et $(\gC',X)$ de matrices respectives $\EE$ et $\EE'$  sont associ\'es, les groupes $G(\EE)$ et $G(\EE')$ sont \'egaux.
\end{propf}
\Dem. 
Le point 1 est imm\'ediat. Les points 2 et 3 se montrent
facilement en remarquant que, pour $\gs$ dans $\SSS_X$, les matrices
$\EE$ et $^\gs \EE$ sont associ\'ees \si il existe une suite $(\nu_1,\ldots,\nu_n)$ de
 coefficients, tous pris dans $\{-1,1\}$ tels que  \\[0.2cm]
$(3)$ \vspace{0,2cm}\centerline{\hspace{3cm}$\forall (i, j)\in X\x X ,
  \qquad \gve_{\gs(i),\gs(j)}=\nu_i\nu_j.\gve_{i,j}$ \hspace{3cm}
  $\Box$}
{\it Notation } :  Lorsque le graphe $\gC$ ou sa matrice $\EE$ sont
d\'etermin\'es  sans ambigu\ii t\'e on notera plus simplement $G$ le
groupe $G(\EE)=G(\gC)$.

\subsection{Localisation} 

Certains graphes associ\'es \`a $(\gC,X)$ vont
jouer un r\^ole important.

\begin{propf}[et d\'efinition de la localisation]$\ $\\ 
Soit $(\gC,X)$ un graphe de matrice $\EE$, et $j$ un point de $X$.\\
$1$. Il existe une  unique matrice $\EE'=(\gve'_{k,l})$ associ\'ee \`a $\EE$ telle que pour tout $k$ dans $X$, $\gve'_{k,j}=1$.  On a\\[0.2cm]
 \vspace{0,2cm}\centerline{ $\forall k,l \in X$, \quad  $\gve'_{k,l}=\nu_k\nu_l.\gve_{k,l}$ \ o\`u \ $\nu_k=\gve_{k,j}$ \ si \ $k\neq j$ \ et \ $\nu_j=1$.}
$2$. Il existe un unique graphe $\gC'$ associ\'e \`a $\gC$ tel que  $j$ soit un point isol\'e de $\gC'$, \cad qu'aucune ar\`ete de $\gC'$ ne contient $j$.\\
Nous dirons que $\EE'$ est la matrice {\rm localis\'ee} de $\EE$ en $j$ et la noterons $\EE'=  {^j\EE}$.\\
De \mm $\gC'$ est le graphe {\rm localis\'e}  de $\gC$ en $j$ et on le note $\gC'= {^j\gC}$.
\end{propf}

\Dem. Bien entendu la deuxi\`eme affirmation n'est que la traduction
de la premi\`ere dans le langage des graphes. Montrons donc la
premi\`ere. Posons pour tout indice $k$ dans $X$,\\[0,1cm]
  \vspace{0,2cm}\centerline{ $\nu_k=\gve_{j,k} $ \  si \ 
$k\neq j$ \ et \   $\nu_j=1$,} 
puis pour deux indices arbitraires $k,l$ dans $X$, \\[0,1cm]
  \vspace{0,2cm}\centerline{$\gve'_{k,l}=\nu_k.\nu_l.\gve_{k,l}$.}
 La matrice $\EE'=(\gve'_{k,l})$ est associ\'ee \`a $\EE$ et satisfait \`a \  \vspace{0.1cm}$\gve'_{j,k}=\nu_j.\nu_k.\gve_{j,k}=\gve_{j,k}^2=1$, pour tout
$k$ dans $X$ ;
elle r\'epond donc \`a la question.  Si $\EE''=(\gve''_{k,l})$ est une deuxi\`eme
solution de coefficient g\'en\'eral
$\gve''_{k,l}=\mu_k.\mu_l.\gve_{k,l}$, o\`u les $\mu_k$ sont pris dans
$\{-1,1\}$, il vient,  pour deux indices $k,l$ dans $X$,\\
\begin{tabular}{ll}
 &$ \gve_{k,l}=\nu_k.\nu_l.\gve'_{k,l}=\mu_k.\mu_l.\gve''_{k,l}$
   \quad et \quad  $ \gve_{k,j}=\nu_k.\nu_j=\mu_k.\mu_j$,\\[0,1cm]
donc & $\gve_{k,j}. \gve_{l,j}=\nu_k.\nu_j.\nu_l.\nu_j=\mu_k.\mu_j.\mu_l.\mu_j$
 \quad  puis \quad  $\nu_k.\nu_l=\mu_k.\mu_l$,\\[0,1cm]
d'o\`u l'on tire & $\gve'_{k,l}=\gve''_{k,l}$ \quad  et enfin \quad  $\EE'=\EE''$.\hspace{5cm}  $\Box$
\end{tabular}

Rassemblons quelques propri\'et\'es de la localisation.

\begin{propf}$\ $\\ 
Soient $\gC$ et $\gC'$ deux graphes sur $X$ de matrices respectives
$\EE$ et $\EE'$.\\
$1$. Pour deux indices $j,k$ dans $X$, on a  \  ${^k\EE}= {^k(}{^j\EE})$, que l'on note plus simplement \ ${^{kj}\EE}$.\\
$2$. Les  matrices $\EE$ et $\EE'$ sont associ\'ees  \si il existe un indice  $k$ dans $ X$ tel que  $^k\EE=  {^k\EE'}$ et dans ce cas l'\'egalit\'e \ $^j\EE=  {^j\EE'}$ \  a lieu pour tous les  indices $j$ de $X$.\\
$3$. Pour toute permutation  $\gs $ de $X$, on a les \'equivalences
entre :\\
$(a)$ Les matrices  $\EE$ et ${^\gs\EE}$ sont associ\'ees.\\
$(b)$ Pour  tout indice  $j$ dans $X$ on a \ $^\gs(^j\EE)=  {^{\gs(j)}\EE}$.\\
$(c)$ Il existe un indice $k$ dans $X$ tel que $^\gs(^k\EE)=  {^{\gs(k)}\EE}$.
\end{propf}
\Dem\\
$1$. Les matrices $\EE$, ${^k\EE}$ et ${^k(}{^j\EE})$ sont associ\'ees
et la $k$-\`eme colonne de ${^k(}{^j\EE})$ ne
contient que le nombre $1$. Or, d'apr\`es la proposition $2$, cette propri\'et\'e caract\'erise la
matrice ${^k\EE}$. Donc ${^k\EE}= {^k(}{^j\EE})$.\\
$2$. Si  $\EE$ et   $\EE'$ sont associ\'ees, alors ${^k\EE}$ et
${^k\EE'}$ le sont aussi mais de plus, leur \mbox{ $k$-\`eme} colonne ne contenant que
des $1$, elles sont \'egales d'apr\`es la \mbox{proposition 2.} Inversement
si \  $^k\EE=  {^k\EE'}$ alors par transitivit\'e  on
voit que $\EE'$ est associ\'ee \`a $\EE$. 
 Dans ce cas on obtient pour tout indice $j$ dans $X$, \\[0.2cm]
 \vspace{0,2cm}\centerline{ $^j\EE= {^{jk}\EE} ={^{jk}\EE'}= {^j\EE'}$.}
$3$. Tout d'abord la d\'efinition 1 nous montre que,  pour
toute permutation $\gs$ de $X$, si les matrices $\EE$ et $\EE'$ sont
associ\'ees alors ${^\gs\EE}$ et ${^\gs\EE'}$ le sont  aussi.\\
Partant de l'hypoth\`ese $(a)$, comme pour chaque indice $j$ dans $X$,
$\EE$ et  $^j\EE$ sont associ\'ees, on en d\'eduit que $\EE$,
${^\gs\EE}$ et $^\gs(^j\EE)$ sont associ\'ees, or \\[0.2cm]
 \vspace{0,2cm}\centerline{$\forall i\in X,\quad
   ^\gs(^j\EE)_{\gs(j),i}={^j\EE}_{\gs^{-1}\gs(j),\gs^{-1}(i)}={^j\EE}_{j,\gs^{-1}(i)}=1$,}
donc $^\gs(^j\EE)$ est l'unique matrice associ\'ee \`a $\EE$ dont la
$\gs(j)$-\`eme colonne ne contient que des $1$, \cad  $^\gs(^j\EE)=
{^{\gs(j)}\EE}$, ce qui prouve $(b)$. \\
Bien entendu $(b)$ implique $(c)$.\\
Supposons enfin que, pour un indice $k$, on ait $^\gs(^k\EE)=
{^{\gs(k)}\EE}$. Alors comme $\EE$,  ${^k\EE}$ et ${^{\gs(k)}\EE}$
sont associ\'ees,   $^\gs\EE$  est associ\'ee \`a
$^\gs(^k\EE)={^{\gs(k)}\EE}$, et par transitivit\'e  $\EE$ et
$^\gs\EE$ sont associ\'ees. Donc $(c)$ implique $(a)$.\hfill $\Box$

{\it Remarque} : en utilisant la correspondance naturelle entre les graphes et leurs matrices on obtient une adaptation imm\'ediate de cette proposition en terme de graphes que nous laissons au lecteur.

Nous allons maintenant interpr\'eter g\'eom\'etriquement l'incidence
de la localisation sur un graphe, c'est \`a dire comparer, pour deux
points $x$ et $y$ de $X$, les  graphes localis\'es ${^x\gC}$ et ${^y\gC}$.\\
 {\it Introduisons une notation }: si $x$ est un point d'un graphe $(\gC,X)$ et $d$ un entier positif, on note $\gC(x,d)$ l'ensemble des points de $X$ qui sont \`a distance $d$ de $x$, et  $\gC(x,d^+)$ l'ensemble de ceux qui sont \`a une distance sup\'erieure ou \'egale \`a $d$.
 
\begin{lemf}$\ $\\ 
Soient $x$ et $y$ deux sommets distincts d'un graphe  $(\gC,X)$ de matrice $\EE$. Notons ${^x\EE}$ et \ ${^y\EE}$ les matrices  des localis\'es ${^x\gC}$ et \ ${^y\gC}$ du graphe $\gC$ en $x$ et $y$.\\
$1.a.$ Il existe une suite de coefficients $\nu_k$ ($k\in X)$, tous dans $\{-1,1\}$, tels que \\[0.2cm]
\vspace{0,2cm}\centerline{ $\forall k,l \in X,  \  {^x\EE}_{k,l}=\nu_k\nu_l.{^y\EE}_{k,l}$, \  o\`u \ $\nu_k=-1 \ \Lra \ k\in {^y\gC}(x,1)\  \Lra \ k \in {^x\gC}(y,1)$}
$1.b.$ On a ${^y\gC}(x,1)={^x\gC}(y,1)$ et les graphes
induits par \  ${^x\gC}$ et \  ${^y\gC}$ sur cet ensemble sont
\'egaux.\\ 
 $2$. Soit $\{k,l\}$ une paire de points dans $X$. La  nature de la liaison entre $k$ et $l$ (ar\`ete ou non) est identique dans les graphes $^x\gC$ et \ $^y\gC$ \si le cardinal de $\{k,l\}\cap {^y\gC}(x,1)$ est pair.
\end{lemf}
\Dem\\
1a. Comme  la matrice ${^x\EE}$ est la localis\'ee en $x$ de ${^y\EE}$ (car ${^x\EE}={^{xy}\EE}$ d'apr\`es la proposition $3$), la proposition $2$ nous montre que les matrices  ${^x\EE}$ et ${^y\EE}$ sont li\'ees par les relations\\[0.1cm]
\vspace{0,2cm}\centerline{ $\forall k,l \in X,  \quad  {^x\EE}_{k,l}=\nu_k\nu_l.{^y\EE}_{k,l}$, \  o\`u \qquad $\nu_k= {^y\EE}_{k,x}$ \  si \ $k\neq x$ \ et \ $\nu_x=1$}
On en d\'eduit que $\nu_k=-1$ \si $\{k,x\}$ est une ar\`ete de ${^y\gC}$, autrement dit si $k\in {^y\gC(x,1)}$.  Or l'\'egalit\'e ${^x\EE}_{k,l}=\nu_k\nu_l.{^y\EE}_{k,l}$ est sym\'etrique en $x$ et $y$ (${^y\EE}_{k,l}=\nu_k\nu_l.{^x\EE}_{k,l}$ est obtenue en multipliant ses deux membres par $\nu_k\nu_l$), donc $\nu_k= {^y\EE}_{k,x}={^x\EE}_{k,y}$, ce qui montre que $\{k,x\}$ est une ar\`ete de ${^y\gC}$ \si $\{k,y\}$ est une ar\`ete de ${^x\gC}$, autrement dit ${^y\gC}(x,1)={^x\gC}(y,1)$,
et ceci prouve 1.a. \\
Revenons \`a l'interpr\'etation g\'eom\'etrique de l'\'egalit\'e $ {^x\EE}_{k,l}=\nu_k\nu_l.{^y\EE}_{k,l}$. Elle signifie 
que la liaison entre les points ${k}$ et $l$ (ar\`ete ou non-ar\`ete)
est de \mm nature dans  les graphes  localis\'es ${^x\gC}$ et
${^y\gC}$ \si le produit $\nu_k\nu_l$ vaut $1$.  On en d\'eduit
imm\'ediatement 1.b. en choisissant $k$ et $l$ dans
${^y\gC}(x,1)={^x\gC}(y,1)$ et plus g\'en\'eralement $2$, en prenant
pour $\{k,l\}$ une paire arbitraire de points dans $X$.\hfill $\Box$
 \begin{center}
  \includegraphics[scale=0.30]{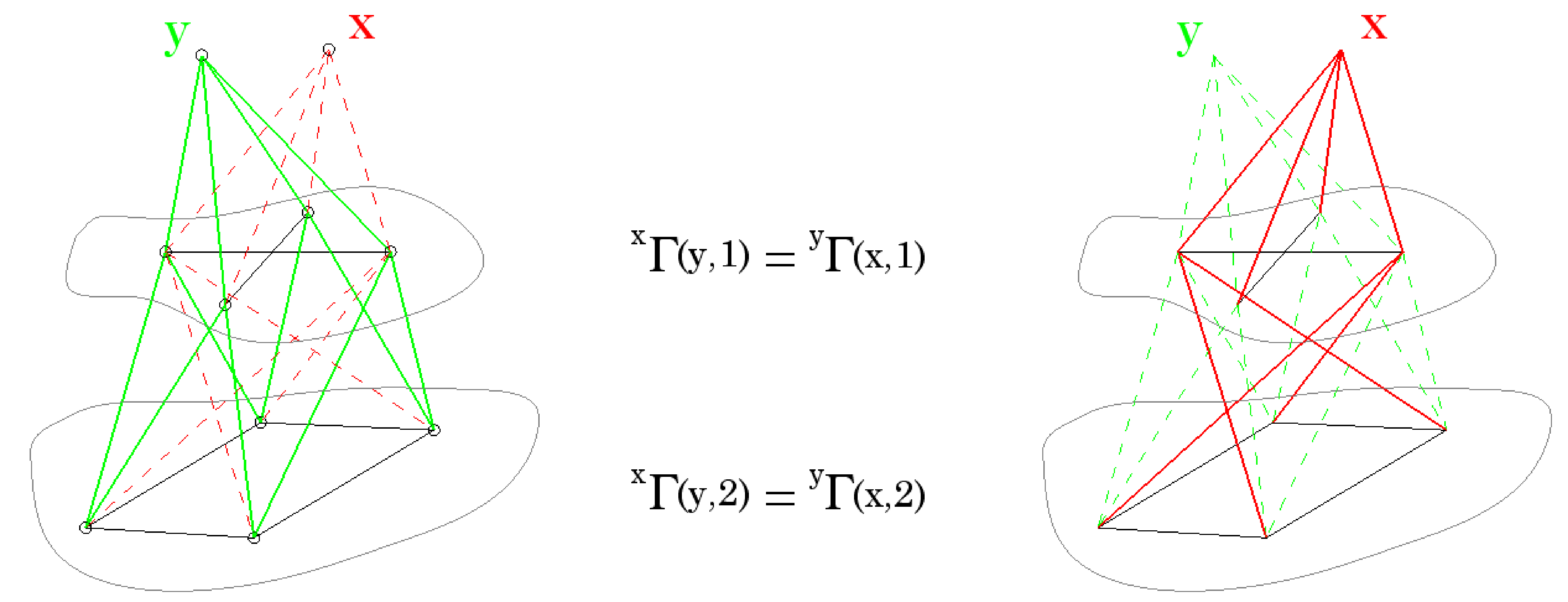}
\npb \centerline{Les graphes $^x\gC$ et $^y\gC$}
 \end{center}
 \subsection{Le probl\`eme de la transitivit\'e} 
 Un graphe $(\gC,X)$ de matrice $\EE$ \'etant donn\'e, la proposition suivante pr\'ecise \`a quelles conditions le groupe $G=G(\EE)$ op\`ere transitivement ou doublement transitivement sur $X$.

\begin{propf}$\ $\\ 
$1$. Le groupe $G$ op\`ere transitivement sur $X$ \si les
localis\'es \ ${^x\gC}$ de $\gC$, pour $x$ variant dans $X$, sont deux \`a deux isomorphes.\\
$2$. S'il op\`ere doublement transitivement sur $X$,  le groupe
d'automorphismes $\Aut({^x\gC})$ de l'un d'eux op\`ere transitivement
sur l'ensemble \ ${X^x}=X-\{x\}$.\\
$3.$ Si pour deux  points  $x$ et $y$ de $X$ distincts  les groupes
d'automorphismes $\Aut({^x\gC})$ et $\Aut({^y\gC})$  op\`erent transitivement, respectivement
sur \ ${X^x}=X-\{x\}$ et \ ${X^y}=X-\{y\}$ alors  le groupe $G$ op\`ere doublement transitivement sur
$X$.
\end{propf}
\Dem. Les points 1 et 2 sont clairs car pour toute permutation
$\gs$ dans le groupe $G$ on a  $^\gs(^x\EE)=  {^{\gs(x)}\EE}$ et
$^\gs(^x\gC)=  {^{\gs(x)}\gC}$ d'apr\`es la proposition
$3$. \\
3. Sous l'hypoth\`ese faite, on voit que
$X-\{x\}$ et $X-\{y\}$ sont respectivement des orbites sur $X$ des
sous-groupes $\Aut({^x\gC})$ et $\Aut({^y\gC})$ de $G$. Mais
l'intersection $(X-\{x\})\cap (X-\{y\})$ n'est pas vide  puisque
$|X|=n\geq 3$, donc le groupe $G$ poss\`ede  une orbite  sur $X$
contenant la r\'eunion $(X-\{x\})\cup (X-\{y\})$ qui n'est autre que
$X$ ; il est donc transitif sur $X$ et comme le sous-groupe  $\Aut({^x\gC})$ de $G$
agit transitivement sur  $X-\{x\}$, on en d\'eduit 
que  $G$ agit deux fois transitivement sur $X$ comme voulu. \hfill $\Box$

Cette proposition nous conduit \`a \'etudier les graphes
satisfaisant aux deux conditions suivantes :\\[0.2cm]
$C1$.  Les localis\'es de $\gC$ en deux points quelconques de $X$ sont isomorphes.\\
$C2$.  Pour un point $x$ de $X$, le groupe d'automorphismes
$\Aut({^x\gC})$ du graphe ${^x\gC}$ op\`ere transitivement sur
l'ensemble \ $X^x=X-\{x\}$.

{\it Notations} : un point $x$ de $X$ \'etant donn\'e, on
note $\gC^x$ le graphe induit par $^x\gC$ sur l'ensemble
${X^x}=X-\{x\}$. Le graphe $\gC^x$ ne diff\`ere donc de $^x\gC$ que
par la suppression du point isol\'e $x$.
 
  \begin{thmf} $\ $\\ 
  Soit $(\gC,X)$ un graphe simple satisfaisant aux conditions $C1$ et $C2$.\\
 $1$. Le groupe $G$ op\`ere deux fois transitivement sur $X$.\\ 
 $2$. Pour tout point $x$ de $X$ le diam\`etre du graphe $(\gC^x,{X^x})$ est au plus $2$, et s'il vaut $1$, $\gC^x$ est le graphe complet sur l'ensemble $X^x$ et $G(\gC)=\SSS_X$.\\
 $3$. Soient $x$ et $y$ deux sommets distincts de $(\gC,X)$\\
 $3.a.$ Les ensembles $\gC^x(y,1)$ et $\gC^y(x,1)$ sont
 \'egaux, de cardinal pair, $2s$.\\
 $3.b$. Les ensembles $\gC^x(y,2)$ et $\gC^y(x,2)$ sont
 \'egaux, de cardinal pair, $2\ov{s}$.\\
 $3.c$. Tout point  de $ \gC^x(y,1)$  est li\'e \`a $\ov{s}$ points
 de $\gC^x(y,2)$ et \`a  $t=2s-\ov{s}-1$ points de $\gC^x(y,1)$.\\
 $3.d$. Tout point de $\gC^x(y,2)$ est li\'e \`a $s$ points de $\gC^x(y,2)$ et \`a $s$ points de  $\gC^x(y,1)$.\\
 $4$. Toute ar\`ete de  $(\gC^x,{X^x})$ est contenue dans exactement $t=2s-\ov{s}-1$ triangles.\\
$5$.  On a $|X|=n=2+2s+2\ov{s}$.
\end{thmf}

\Dem : \\
$1$. Tout d'abord les conditions $C1$ et $C2$ nous montrent que les actions
des groupes $\Aut({^x\gC})$ et $\Aut({^y\gC})$ sur les ensembles
$X^x=X-\{x\}$ et $X^y=X-\{y\}$ sont transitives. Par la
proposition $4$ on en d\'eduit que le groupe $G$ op\`ere deux fois
transitivement sur $X$. Donc les cardinaux des ensembles
$\gC^x(y,1)$, $\gC^x(y,2)$ et $\gC^x(y,2^+)$ ne d\'ependent pas du
choix des sommets $x$ et $y$ pourvu qu'ils soient distincts.\\
$2$. Dire que le diam\`etre du graphe $(\gC^x,{X^x})$ est $1$ \'equivaut \`a dire que c'est le graphe complet sur $X^x$ et   l'\'egalit\'e $G(\gC)=\SSS_X$ en d\'ecoule imm\'ediatement.\\
Comme le graphe $\gC^x$ ne
diff\`ere de $^x\gC$ que par la suppression du point isol\'e $x$
l'\'egalit\'e  $\gC^x(y,1)=\gC^y(x,1)$ vient  simplement du lemme $1\ (1.b)$.\\
 Supposons que le diam\`etre de $(\gC^x,{X^x})$ est au moins
$2$, et choisissons un point $z_2$ dans $\gC^x(y,2^+)$ et 
un point $z$ dans $X$. Comme  $x$ est isol\'e dans $^x\gC$, on a
 $\gC^x(y,1)= {^x\gC}(y,1)$ et le
lemme $1$ nous montre que si $z$ n'est pas dans
$\gC^x(y,1)$ la liaison entre $z$ et $z_2$ est
identique dans les graphes  $^x\gC$ et  $^y\gC$, tandis que si $z$ est dans
$\gC^x(y,1)=\gC^y(x,1)$ (repr\'esent\'e par $z'$ sur la figure A) alors $\{z_2,z\}$ est une ar\`ete de  $\gC^x$
\si c'est une non-ar\`ete de $\gC^y$. Or la double transitivit\'e du
groupe $G$ agissant  sur $X$ montre que le nombre d'ar\`etes issues
de $z_2$ doit \^etre le \mm dans les graphes  $\gC^x$ et  $\gC^y$. On
en d\'eduit donc que le nombre d'\elts de $\gC^x(y,1)$ est un
nombre pair, $2s$, et que $z_2$ est li\'e \`a exactement $s$
d'entre eux dans le graphe $\gC^x$. Mais on en d\'eduit aussi que  $z_2$ est \`a distance $2$
de $y$ dans le graphe $\gC^x$, donc les graphes $\gC^x$ et $\gC^y$ sont de diam\`etre $2$ et les ensembles $\gC^x(y,2)$ et $\gC^y(x,2)$ sont  \'egaux car ils sont le compl\'ementaire commun dans $X$ de  $\{x,y\}\cup \gC^x(y,1)$. De plus les $s$ sommets de $\gC^x$ qui sont li\'es \`a $z_2$ mais pas dans $\gC^x(y,1)$ sont n\'ecessairement dans $\gC^x(y,2)$. Ceci ach\`eve la preuve de   $2$, $3.a$ et $3.d.$
\begin{center}
  \includegraphics[scale=0.25]{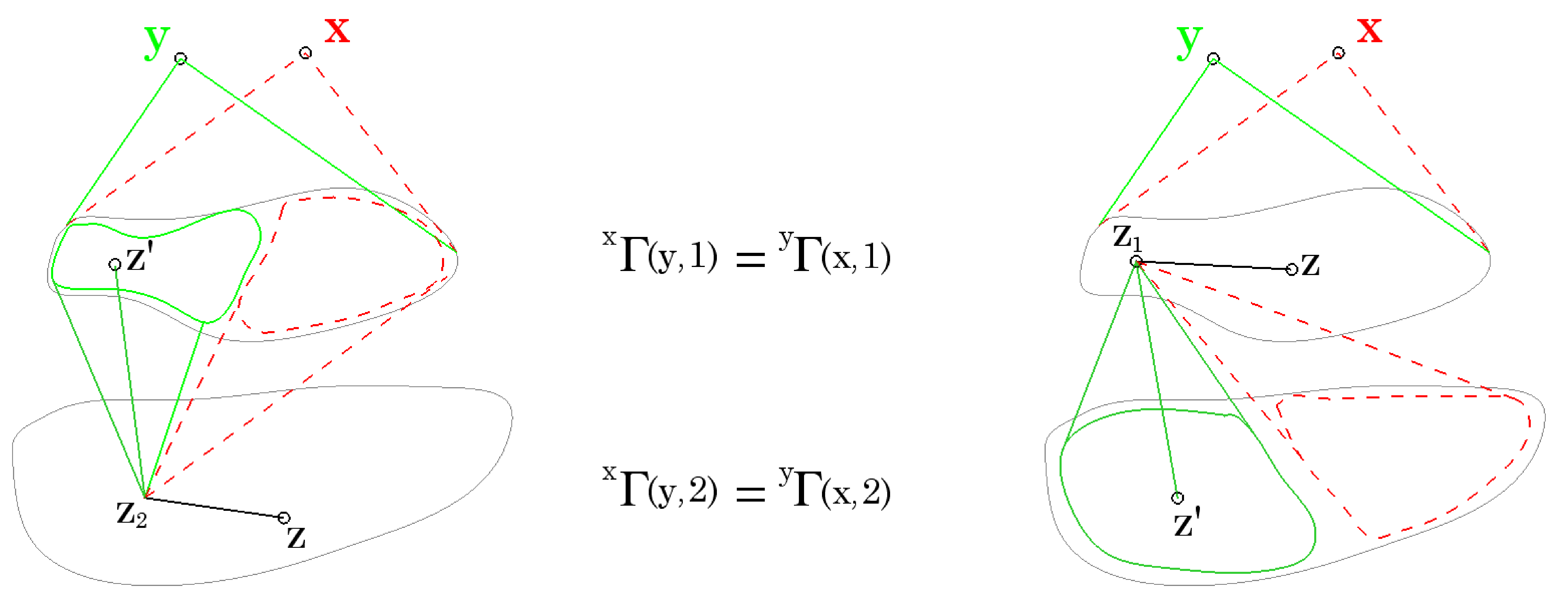}
\npb \centerline{figure A \hspace{8cm} figure B }
 \end{center}
Choisissons maintenant un point $z_1$ dans $\gC^x(y,1)$ et un point
$z$ dans $X$. En raisonnant de fa\cc on analogue, on voit que les
liaisons entre $z_1$ et $z$ sont de \mm nature dans les graphes
$^x\gC$ et $^y\gC$ si $z$ est dans $\gC^x(y,1)$, et de nature contraire
si $z$ est dans $\{x,y\}\cup \gC^x(y,2)$ (repr\'esent\'e par $z'$ sur la  figure B), donc le nombre d'\elts de
$\gC^y(x,2)=\gC^x(y,2)$ est un nombre pair, $2\ov{s}$, et le sommet $z_1$ est li\'e \`a exactement $\ov{s}$
d'entre eux dans dans le graphe $\gC^x$. Tenant compte du fait que $z_1$ est li\'e \`a $y$, on
en d\'eduit, par soustraction que $z_1$ est li\'e \`a exactement
$t=2s-\ov{s}-1$ points de $\gC^x(y,1)$, ce qui ach\`eve la preuve de $3.$
Enfin ce nombre $t=2s-\ov{s}-1$ est aussi le nombre de triangles dont
$\{y,z_1\}$ est une ar\`ete, et ceci prouve $4$.\\
$5$. C'est une simple addition puisque $X=\{x,y\}\cup \gC^x(y,1)\cup \gC^x(y,2)$. \hfill $\Box$

Avant de passer \`a la construction de graphes satisfaisants aux conditions du \ttt 1, prenons un instant pour nous accorder avec les notations usuelles.
Un graphe simple $(\gL,Y)$ est dit {\it fortement  r\'egulier} lorsqu'il poss\`ede les trois propri\'et\'es suivantes :\\
$*$ Chaque sommet de $\gL$ est li\'e \`a exactement $k$ sommets de $\gL$. \\
$*$ Chaque ar\`ete est contenue dans exactement $\gl$ triangles.\\
$*$ Deux points non li\'es sont simultan\'ement li\'es \`a $\mu$ sommets.\\
On note souvent $v$ le cardinal de $Y$ et on dit que le graphe $(\gL,Y)$ est de type $(v,k,\gl,\mu)$. Les param\`etres d'un tel graphe ne sont pas totalement ind\'ependants puisqu'il v\'erifient la relation 
$(v-k-1)\mu = k(k-\gl-1)$. 
Il est facile de v\'erifier que les  conditions impos\'ees par le \ttt 1 au graphe $(\gC^x,X^x)$  en font un graphe fortement r\'egulier avec les param\`etres $ (v,k,\gl,\mu)=(n-1,2s,t,s)$, et l'\'egalit\'e  $(v-k-1)\mu = k(k-\gl-1)$ se traduit alors par $n=2+2s+2\ov{s}$.\\
Pour des raisons essentiellement typographiques (les lettres $v,k,\gl,$ et $\mu$ sont bien utiles) nous gardons dans la suite les param\`etres $n,t,s$ et $\ov{s}$ pour d\'ecrire le graphe $(\gC,X)$.

 \section{Constructions de graphes}
 
 {\it Fixons quelques notations.} \\
Dans cette partie, un point $x$ de $X$ \'etant donn\'e, on renomme $(\gL,Y)$ le graphe $(\gC^x,{X^x})$ et on suppose qu'il satisfait aux conditions $3$, $4$ et $5$ du \ttt 1 :\\
1. Le graphe  $(\gL,Y)$ est diam\`etre $2$.\\
2. Pour tout point $y$ de $Y$, \\
2.a.  L'ensemble $\gL(y,1)$ des points li\'es \`a $y$ est de cardinal pair $2s$.\\
2.b.  L'ensemble $\gL(y,2)$ des points \`a distance $2$ de $y$ est de cardinal pair $2\ov{s}$.\\
2.c.  Tout point $z$ de  $\gL(y,1)$ est li\'e \`a $\ov{s}$ points de $\gL(y,2)$ et \`a  $t=2s-\ov{s}-1$ points de $\gL(y,1)$.\\
 2.d. Tout point $z$ de $\gL(y,2)$ est li\'e \`a $s$ points de $\gL(y,2)$ et \`a $s$ points de  $\gL(y,1)$.\\
 3. Chaque ar\`ete du graphe  $(\gL,Y)$ est contenue dans exactement \ $t=2s-\ov{s}-1$\  triangles.\\
4.  On a $|Y|=1+2s+2\ov{s}$.

Nous dirons qu'un graphe  $(\gL,Y)$ satisfaisant \`a ces quatre propri\'et\'es est  {\it  extensible} et que les constantes $(t,s,\ov{s})$ sont ses {\it param\`etres}. 
Lorsqu'aucune ambigu\ii t\'e n'est possible concernant le point $y$ de $Y$, on notera  $(\gL_1,Y_1)$ et $(\gL_2,Y_2)$ les graphes induits par $(\gL,Y)$ sur les ensembles $Y_1=\gL(y,1)$ et $Y_2=\gL(y,2)$, \cad les graphes dont les ar\`etes sont celles de $(\gL,Y)$ qui sont contenues dans $Y_1$, ou $Y_2$.

Remarquons que la connaissance du graphe $(\gL,Y)=(\gC^x,{X^x})$
d\'etermine, \`a isomorphisme pr\`es, pour tout $y$ dans $X$, le graphe $(\gC^y,{X^y})$, puisque ${^y\gC}={^{yx}\gC}$. Plus pr\'ecis\'ement on a 

 \begin{propf} $\ $\\  
 Si pour un point $x$ de $X$ le graphe $(\gC^x,{X^x})$ est
 extensible, alors le graphe $(\gC^y,{X^y})$ l'est pour tout $y$ dans $X$, et les param\`etres des graphes $(\gC^x,{X^x})$ et $(\gC^y,{X^y})$ sont \'egaux.
\end{propf}
\Dem. Simple d\'ecompte. \hfill $\Box$

 \subsection{Le graphe compl\'ementaire}
 On appelle  {\it compl\'ementaire} d'un graphe $(\gL,Y)$ le graphe $(\ov{\gL},Y)$  admettant le \mm ensemble de sommets $Y$ que  $(\gL,Y)$ et dont les ar\`etes sont les non-ar\`etes de $(\gL,Y)$. 
 
 \begin{propf} $\ $\\  
Si $(\gL,Y)$ est un graphe extensible, son graphe compl\'ementaire  $(\ov{\gL},Y)$ l'est aussi et ses param\`etres $(\ov{t},s',\ov{s'})$ sont li\'es aux param\`etres $(t,s,\ov{s})$ de  $(\gL,Y)$ par \\[0,1cm]
 \centerline{$s'=\ov{s}$, \quad $\ov{s'}=s$, \quad  $t+\ov{t}=s+\ov{s}-2 =s'+\ov{s'}-2$}
 \end{propf}
 
  Cette proposition, qui se montre tr\`es simplement, nous permet de limiter notre recherche aux cas o\`u $t\leq \ov{t}$, ce qui \'equivaut \`a \ $2t\leq s+\ov{s}-2$.
  
On aborde maintenant la recherche des graphes extensibles $(\gL,Y)$ de param\`etres $(t,s,\ov{s})$. On organise la discussion en fonction du param\`etre $t$ qui nous donne le nombres de triangles incidents \`a une ar\`ete donn\'ee. Nous regardons tout d'abord les cas $t=0$ et $t=1$. Ensuite on verrons que les graphes de Paley $P(q)$ forme une famille infinie de solutions admettent un syst\`eme de param\`etres de la forme $(t,s,\ov{s})=(t,t+1,t+1)$ o\`u $q=4t+5$. 

 \subsection{$t=0$}
 
  \begin{thmf} $\ $\\ 
$1.$ Si $t=0$, le graphe $(\gL,Y)$  existe uniquement si $(t,s,\ov{s})=(0,1,1)$. C'est un pentagone.\\
$2.$ Son groupe d'automorphismes $\Aut(\gL)$ est le groupe di\'edral $D_5$, 
dont l'extension $G(\gC)$ est un groupe op\'erant doublement
transitivement sur l'ensemble $X=Y\cup\{x\}$, semblable \`a la repr\'esentation doublement transitive du groupe altern\'e $A_5$ sur un ensemble de cardinal $6$.
\end{thmf}
\Dem \\
1. Il est clair qu'un pentagone est un graphe extensible de param\`etre $t=0$.\\
 Inversement si $(\gL,Y)$ est extensible et $t=0$, choisissons un
 point $y$ dans $Y$ et un point  $z$ dans $Y_2=\gL(y,2)$ qui est donc
 li\'e \`a $2s$ points de $Y$ dont $s$ sont dans $Y_1=\gL(y,1)$ et les
 autres dans $Y_2$. Comme $(\gL,Y)$ ne contient aucun triangle chaque
 point $z_2$ de $\gL(z,1)\cap Y_2$ est li\'e \`a $s$  points dans
 $Y_1$ (propri\'et\'e 2d des graphes extensibles) qui sont tous dans le compl\'ementaire $C$ de $\gL(z,1)\cap Y_1$ dans $Y_1$. Or, d'apr\`es la condition 2d, $|\gL(z,1)\cap Y_1|=|C|=s$, donc inversement chaque chaque point de $C$ est li\'e aux $s$ points de $\gL(z,1)\cap Y_2$. \\
Supposons que $s>1$. On peut alors choisir deux points distincts $y_1$
et $z_1$ dans $C$ li\'es \`a $y$ et aux $s$ points de $\gL(z,1)\cap
Y_2$. Il vient donc \mbox{$|\gL(y_1,1)\cap \gL(z_1,1)|\geq s+1$}, ce qui contredit la condition 2d. Donc, par l'absurde, $s=1$, $\ov{s}=1$ et le graphe $(\gL,Y)$ est un pentagone.  \\
2. Revenons au graphe $(\gC,X)$ pour lequel
$(\gC^x,X^x)=(\gL,Y)$. D'apr\`es la proposition 5, pour tout $y$ dans
$X$, le graphe $(\gC^y,X^y)$ est extensible, de m\^emes param\`etres que  $(\gC^x,X^x)$, et l'unicit\'e de la construction ci-dessus nous montre qu'il s'agit aussi d'un pentagone. Le groupe $G(\gC)$ op\`ere donc transitivement sur $X$ et le stabilisateur d'un point $x$ de $X$ est semblable, dans son action sur $X^x$ au groupe di\'edral $D_5$.
 \hfill $\Box$ 

 \subsection{$t=1$}
 
 Nous allons montrer que lorsque le param\`etre  $t$  vaut $1$, il existe, \`a isomorphisme pr\`es,  quatre graphes  extensibles $(\gL,Y)$ dont les param\`etres sont \\[0,1cm]
 \vspace{0,1cm}\centerline{$ (1,1,0)$, \quad $ (1,2,2)$, \quad $ (1,3,4)$, \quad $ (1,5,8)$}
Comme pr\'ec\'edemment on choisit un sommet  $y$ dans $Y$ et on note $(\gL_1,Y_1)$ et $(\gL_2,Y_2)$, les graphes induits par $(\gL,Y)$ sur  $Y_1=\gL(y,1)$ et $Y_2=\gL(y,2)$.
 
   \begin{propf} $\ $
 Pour des indices $i$ et $j$  entre $1$ et $s$ on a :\\
 $(1)$ Le graphe $(\gL_1,Y_1)$  dessine sur  $Y_1$ une  partition  par $s$ ar\`etes, qu'on notera \\[0,1cm]
   \vspace{0,2cm}\centerline{$\ga_1=\{a'_1,a''_1\}$,$\ldots$,$\ga_s=\{a'_s,a''_s\}$.}
 $(2)$ Notant  $A'_i$ (resp. $A''_i$) l'ensemble des points de $Y_2$  li\'es \`a $a'_i$ (resp. $a''_i$), pour chaque indice $i$,  $\{A'_i,A''_i\}$ est une partition de
  $Y_2$ en deux sous-ensembles de cardinal \ $\ov{s}=2(s-1)$. \\
  $(3)$  Le graphe $(\gL_2,Y_2)$ ne contient aucun triangle.\\
 $(4)$ Toute ar\`ete de $(\gL_2,Y_2)$ est contenue dans exactement un ensemble parmi $A'_1, A''_1,\ldots, A'_s, A''_s$, et le graphe  $(\gL_2,Y_2)$ induit  sur chacun d'eux une partition par $s-1$ ar\`etes.\\
  $(5)$  Si  $\ga'$ et $\ga''$ sont deux ar\`etes du graphe $(\gL_2,Y_2)$ choisies  respectivement dans  $A'_i$ et $ A''_i$, il existe des ar\`etes $\gb'$ et $\gb''$ du graphe $(\gL_2,Y_2)$, uniquement d\'etermin\'ees telles que $\ga'$, $\ga''$, $\gb'$ et $\gb''$ forment un carr\'e. \\
  De plus  il existe un unique indice $j$ tel que $\gb'$ et $\gb''$ soient  respectivement contenues dans  $A'_j$ et $ A''_j$.
 \end{propf}
 \begin{center}
  \includegraphics[scale=0.25]{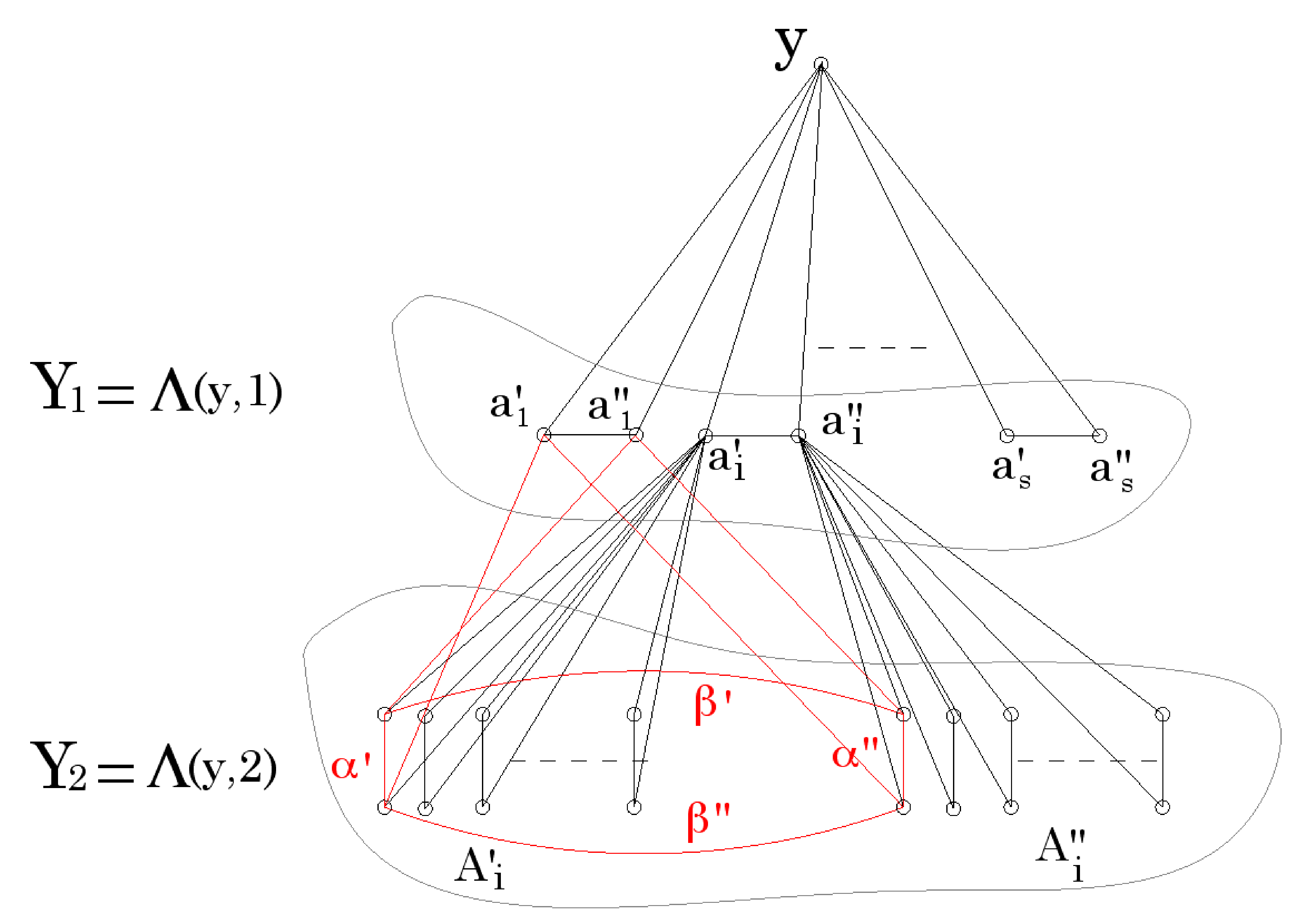}
\npb \centerline{Repr\'esentation symbolique de la proposition 7}
 \end{center}
\Dem \\ 
(1) Comme $t=1$, l'ensemble des triangles de sommet $y$ dessine sur $Y_1$ une partition par $s$ ar\`etes, que l'on note $\ga_1=\{a'_1,a''_1\},\ldots, \ga_s=\{a'_s,a''_s\}$.\\
(2)  Chacun des sommets $a'_i$ (resp. $a''_i$) est donc li\'e \`a un ensemble $A'_i$ (resp. $A''_i$) de $2s-2$ points dans $Y_2$. La condition $t=1$ nous montre  que pour chaque indice $i$ l'intersection  $A'_i \cap A''_i$ est vide et donc  $\{A'_i,A''_i\}$ est une partition de $Y_2$ \\
 (3) S'il existait un triangle  $T$ dans $Y_2$, une de ses ar\`etes serait contenue dans l'un des ensembles $A'_1$ ou $ A''_1$ et serait donc commune \`a deux triangles ($T$ et un triangle de sommet $a'_1$ ou $a''_1$), ce qui contredirait la condition $t=1$. Ceci prouve $(3)$.\\
(4) La condition  $t=1$  nous montre aussi que, pour chaque indice  $i$, l'ensemble des traces sur $Y_2$ des triangles de sommet $a'_i$ (resp. $a''_i$)  dessine sur  $A'_i$ (resp. $A''_i$)
 une partition par $s-1$ ar\`etes. Inversement toute ar\`ete contenue dans $Y_2$ est contenue dans un unique triangle dont le troisi\`eme sommet est n\'ecessairement l'un des $a'_i$ ou des $a''_i$. Donc cette ar\`ete est contenue dans  $A'_i$ ou $ A''_i$. \\
 (5) Soient $u$ et $v$ les deux sommets de l'ar\`ete $\ga'$ et
 $\beta_j=\{u,u_j\}$  la liste des $s-1$ ar\`etes d'origine $u$
 ($1\leq j \leq s-1$), contenues dans $Y_2$ et distinctes de
 $\ga'$. Comme, d'apr\`es $(4)$, aucune d'elle ne peut \^etre contenue
 dans $A'_i$, leurs extr\'emit\'es $u_1,\ldots, u_{s-1}$ sont toutes
 contenues dans $ A''_i$. Mais comme $Y_2$ ne contient aucun triangle
 l'ensemble $\{u_1,\ldots, u_{s-1}\}$ ne contient aucune ar\`ete.
 Donc chacune des $s-1$ ar\`etes contenue dans $ A''_i$ contient
 exactement l'un des points $u_j$. En particulier il existe un unique
 indice $k$ tel que $u_k$ soit un sommet de l'ar\`ete $\ga''$. En raisonnant de \mm avec
 le point $v$, on voit qu'il existe un sommet $v_l$ de l'ar\`ete
 $\ga''$ tel que $\{v,v_l\}$ soit une ar\`ete. De plus comme $Y_2$ ne
 contient pas de triangle, on voit que $u_k\neq v_l$, autrement dit
 $\ga''=\{u_k,v_l\}$. Ceci nous prouve la premi\`ere partie de $(5)$
 (avec $\gb'=\{u,u_k\}$ et $\gb''=\{v,v_k\}$). Soient maintenant $a_k$
 (resp. $a_l$) le troisi\`eme sommet du triangle dont  $\gb'$
 (resp. $\gb''$) est une ar\`ete. Comme les deux sommets $u,v$ de ces
 triangles sont li\'es par une ar\`ete, le raisonnement qui
 pr\'ec\`ede nous montre que le graphe induit sur les quatre sommets
 $u_k,v_l,a_k,a_l$ est un carr\'e. Or $\{a_k,u_k\}$, $\{u_k,v_l\}$ et
 $\{v_l,a_l\}$ sont des ar\`etes donc $\{a_k,a_l\}$ est la quatri\`me ar\`ete,
 ce qui montre que pour un indice $j$ bien d\'etermin\'e on a
 $(a_k,a_l)=(a'_j,a''_j)$ et ceci termine la preuve de $(5)$.\hfill $\Box$ 

{\it Remarques et terminologie.}  Lorsque $s=1$, on a $\ov{s}=2(s-1)=0$ et l'\'etude du graphe $(\gL_2,Y_2)$ perd tout int\'er\^et ; on supose donc maintenant que $s\geq 2$. D'apr\`es la propri\'et\'e  $(4)$,  \`a chaque ar\`ete $\ga$ du graphe $(\gL_2,Y_2)$ on peut associer un unique indice  $i$, $(1\leq i\leq s)$, tel que $\ga$ soit contenue dans $A'_i$ ou dans $A''_i$. L'ensemble des ar\`etes ainsi associ\'ees \`a un \mm indice forment une partition $D_i$ de $Y_2$ qu'on appelle leur {\it direction}. \`A chaque direction $D_i$ on associe aussi la {\it translation } de $Y_2$ donn\'ee par $t_i(u)=v$ \si $\{u,v\}$ est une ar\`ete de direction $D_i$. Enfin nous dirons qu'une partie $S$ de $Y_2$ est {\it satur\'ee} si elle est stable par toute translation dont la direction est celle d'une ar\`ete contenue dans $S$. Ainsi, par exemple, toute ar\`ete est une partie satur\'ee et la proposition 7 nous montre que les carr\'es sont des parties satur\'ees de  $(\gL_2,Y_2)$.  En reprenant ses notations  on peut  pr\'eciser la structure du graphe $(\gL_2, Y_2)$ :

\begin{lemf}  Choisissons deux indices $i$ et $j$  entre $1$ et $s$. \\
$(a)$ Le groupe $T$ des permutations de $Y_2$ engendr\'e par les translations $t_1,\ldots,t_s$ est ab\'elien $2$-\'el\'ementaire, et op\`ere r\'eguli\`erement sur $Y_2$.\\
$(b)$ Il est de rang  $r\geq s-1$, et s'il est de rang $s-1$ alors $t_s=t_1\oo \ldots \oo t_{s-1}.$\\
$(c) $ On a $(t,s, \ov{s})\in \{(1,1,0),  (1,2,2), (1,3,4),  (1,5,8)\}$ 
\end{lemf}

  \Dem \\
  $(a)$ Consid\'erons  $t_i$ et $t_j$ deux translations de directions $D_i$ et $D_j$ distinctes, $u$ un point de $Y_2$ et    les ar\`etes\\[0,1cm]
  \vspace{0,2cm}\centerline{  $\ga'_i=\{u,t_i(u)\}$, $\ga''_i=\{t_j(u),t_it_j(u)\}$,  $\gb'_j=\{u,t_j(u)\}$, $\gb''_j=\{t_i(u),t_jt_i(u)\}$.}
   Les ar\`etes $\ga'_i$ et $\ga''_i$ ont m\^eme direction, mais ne sont pas simultan\'ement dans $A'_i$ ou dans $A''_i$, sinon l'ar\`ete  $\gb'_j$ serait elle \mm contenue dans $A'_i$ ou $A''_i$, ce qui contredirait $(4)$. La propri\'et\'e $(5)$ nous montre donc l'existence d'une unique ar\`ete $\gb''$  de \mm direction que $\gb'_j$ telle  que les quatre ar\`etes $\ga'_i$, $\ga''_i$, $\gb'_j$ et  $\gb''$ forment un carr\'e. Or  ceci implique clairement  que  $\gb''_j=\gb''$, puis que $t_i\oo t_j(u)= t_j\oo t_i(u)$. Cet argument s'appliquant \`a tout point $u$ de $Y_2$, on en d\'eduit que  $t_i\oo t_j= t_j\oo t_i$. Le groupe $T$, qui est donc commutatif et engendr\'e par des involutions, est  $2$-\'el\'ementaire. Enfin la condition $(5)$ de la proposition 7 nous montre que $T$ op\`ere transitivement sur $Y_2$, mais aussi r\'eguli\`erement puisqu'il est commutatif.\\
  $(b)$ Soit $r$ le rang de $T$ et, quitte \`a renommer les translations $t_1, \ldots, t_s$, on suppose que $t_1, \ldots, t_r$ est un syst\`eme g\'en\'erateur minimal de $T$.  On se donne un point $u$ dans $Y_2$ et \`a chaque ar\`ete $\ga_i$ d'origine $u$ et de direction $D_i$ on associe le coefficient $\gve_i=\gve_i(u)$ qui vaut $1$ si l'ar\`ete $\ga_i$ est contenue dans $A'_i$ et $-1$ si elle est dans $A''_i$. On note $\gve(u)=(\gve_1,\ldots,\gve_s)$. \\
  Compte tenu de la description des carr\'es donn\'ee dans la proposition 8, $(5)$, on v\'erifie que la translation de vecteur $t_j$, $(1\leq j\leq s)$, transforme le point $u$ en une image $t_j(u)$ associ\'ee au $s$-uplet  $\gve(t_j(u))=(-\gve_1,\ldots,-\gve_{j-1},\gve_j, -\gve_{j+1}\ldots,-\gve_s)$. Choisissons   un indice $j>r$ et supposons, quitte \`a modifier l'ordre des g\'en\'erateurs $t_1,\ldots, t_r$ de $T$, que pour un entier $k$ compris entre $2$ et $r$ on ait  $t_j=t_1\oo\ldots\oo t_k$. En \'evaluant de deux fa\cc ons $\gve(t_j(u))$ il vient \\[0,1cm]
     \vspace{0,2cm}\centerline{$(-\gve_1,\ldots,-\gve_{j-1},\gve_j, -\gve_{j+1}\ldots,-\gve_s)= ((-1)^{k-1}\gve_1,\ldots,(-1)^{k-1}\gve_{k},(-1)^{k}\gve_{k+1},\ldots,(-1)^{k}\gve_s)$}
 Mais alors en utilisant les in\'egalit\'es $2\leq k\leq r<j\leq s$, on voit  successivement que \ 
 $k$ est pair, $j=s$ (donc la deni\`ere \'egalit\'e '' impossible ''  $-\gve_s=(-1)^k\gve_s$ n'a pas lieu), $s=r+1$ puisque $j>r$ implique $j=s$ et enfin $k=r$, ce qui prouve $(2)$.\\
 $(c)$ Le groupe $T$, de rang $r$,  op\'erant r\'eguli\`erement sur $Y_2$ qui est de cardinal $2.\ov{s}=4(s-1)$, on a $2^r=2.\ov{s}=4.(s-1)$ et $s-1\leq r \leq s$ \ d'apr\`es ce qui pr\'ec\`ede.  Donc $s$ doit v\'erifier l'une des deux \'equations $2^s=4.(s-1)$ ou $2^{s-1}=4.(s-1)$. La premi\`ere conduit aux solutions $ (1,2,2)$ et  $(1,3,4)$, la deuxi\`eme conduit \`a la solution $(1,5,8)$. Enfin la solution $(1,1,0)$ est particuli\`ere puisqu'elle correspond au cas o\`u $Y_2$ est vide ; dans ce cas le graphe $(\gL,Y)$ est un triangle.\nlb\hfill$\Box$
 
 Venons en \`a la synth\`ese : le \ttt 2  nous montre que les renseignements obenus sur le graphe $(\gL_2,Y_2)$ dans le lemme 2 sont suffisant pour \'etablir l'existence et l'unicit\'e des graphes ayant ces propri\'et\'es.
 
 \begin{thmf} $\ $ \\ 
$1.$  Lorsque $t=1$, il existe, \`a isomorphisme pr\`es, quatre
graphes $(\gL,Y)$ extensibles.  Leurs param\`etres $(t,s,\ov{s})$ sont :\\[0,1cm]
\vspace{0,2cm}\centerline{$ (1,1,0)$, \quad $ (1,2,2)$, \quad $ (1,3,4)$, \quad $ (1,5,8)$}
$2.a.$  Si $(t,s,\ov{s})=(1,1,0)$, on a $|Y|=1+2s+2\ov{s}=13$,  et $|X|=|Y|+1=4$. \\
Le graphe $(\gL,Y)$ est un triangle. Le groupe $\Aut(\gL,Y)$ est isomorphe \`a $\SSS_3$ et le groupe $G(\gC)$ est isomorphe \`a $\SSS_4$.\\
$2.b.$    Si $(t,s,\ov{s})=(1,2,2)$, on a  $|Y|=9$  et $|X|=|Y|+1=10$.  \\
Le groupe $\Aut(\gL_2,Y_2)$ est semblable au groupe $C$ des isom\'etries d'un carr\'e, donc d'ordre $|C|=8$,  et le groupe $\Aut(\gL,Y)$ \'etant transitif sur $Y$ il est d'ordre $9.|C|=72$. Enfin le groupe $G(\gC)$, doublement transitif sur $X$, est d'ordre $10.72=720$.\\
$2.c.$    Si $(t,s,\ov{s})=(1,3,4)$, on a $|Y|=15$ et $|X|=|Y|+1=16$.\\
 Le groupe $\Aut(\gL_2,Y_2)$ est semblable au groupe $C$ des isom\'etries d'un cube, donc  d'ordre $|C|=48$,  et le groupe $\Aut(\gL,Y)$ \'etant transitif sur $Y$, il est d'ordre $15.|C|=720$. Enfin le groupe $G(\gC)$, doublement transitif sur $X$,  est d'ordre $16.720=11520$.\\
$2.d.$    Si $(t,s,\ov{s})=(1,5,8)$,    on a $|Y|=27$ et $|X|=|Y|+1=28$.\\
 Le groupe $\Aut(\gL_2,Y_2)$ est semblable au groupe $H$ des
 isom\'etries d'un hypercube en dimension $4$ compl\'et\'e par ses
 diagonales principales, donc d'ordre $|H|=2^4.5!=1920$,  et le groupe
 $\Aut(\gL,Y)$ \'etant transitif sur $Y$, il est d'ordre
 $27.|H|=51840$. Enfin le groupe $G(\gC)$, doublement transitif  sur
 $X$, est d'ordre $28.51840=1451520$. Il est
 semblable \` a la repr\'esentation de degr\'e $28$ du groupe simple $S_6(2)=PSp_6(2)$.
\end{thmf}
\Dem \\ 
Le principe est assez simple. L'analyse pr\'ec\'edente nous a montr\'e que, pour tout point $y$ d'un graphe extensible $(\gL,Y)$ de param\`etres $(t,s,\ov{s})$, le graphe $(\gL_2,Y_2)$ induit par $(\gL,Y)$ sur $Y_2=\gL(y,2)$ satisfait aux conditions de la proposition 7 et du lemme 2. Donc, si inversement, ces conditions nous permettent de construire le graphe $(\gL,Y)$, \`a isomorphisme pr\`es, on en d\'eduira d\'ej\`a l'existence et l'unicit\'e du graphe $(\gL,Y)$ de param\`etres $(t,s,\ov{s})$, mais aussi la transitivit\'e sur $Y$ du groupe de ses automorphismes $\Aut(\gL)$.  En utilisant alors la proposition 5, on voit que les  localis\'es du graphe initial $(\gC,X)$ en des points  distincts de $X$ sont isomorphes et on en d\'eduit la double transiivit\'e de l'action de $G=\Aut(\gC,X)$ sur $X$ en utilisant par exemple la proposition 4.\\
$*$ $s=1$. Ce cas est trivial. On a $Y_2=\o$ et le graphe $(\gL,Y)$ est un triangle. Donc le groupe $\Aut(\gC,X)$ est le groupe $\SSS_4$ agissant naturellement sur $X$.\\
$*$ {\it Autres cas.} Pour montrer l'unicit\'e de ces graphes (\`a isomorphisme pr\`es) on utilise deux remarques :\\
$**$ Le nombre $s$ \'etant donn\'e, ainsi que $\ov{s}=2(s-1)$, dans chaque cas le cardinal de $Y_2$ est une puissance de $2$ et il n'existe qu'une repr\'esentation  r\'eguli\`ere d'un groupe $T$, $2$-\'el\'ementaire sur $Y_2$. Par ailleurs si $t_1,\ldots,t_r$ est un syst\`eme g\'en\'erateur minimal de $T$, le lemme $2$ nous montre que les ar\`etes du graphe $(\gL_2,Y_2)$ sont les orbites sur $Y_2$ de chacune des $t_i$ ($1\leq i\leq r$), et \'eventuellement de la translation   $t_s=t_1\oo \ldots \oo t_{s-1}$ lorsque le rang $r$ de $T$ est $s-1$ (\cad si  $s=5$). Ceci prouve que les liaisons internes au graphe $(\gL_2,Y_2)$ sont enti\`erement d\'etermin\'ees par la connaissance de l'action de $T$ sur $Y_2$.\\
$**$ Pour d\'eterminer compl\`etement le graphe  $(\gL,Y)$ il nous suffit donc de reconstruire les liaisons entre  $(\gL_2,Y_2)$ et  $(\gL_1,Y_1)$ et les liaisons internes \`a  $(\gL_1,Y_1)$. Or les conditions $(1)$ et $(2)$ de la proposition 7, nous montrent que  l'ensemble $\{a'_1,a''_1\ldots, a'_s,a''_s\}$ des points de $Y_1$ est en  bijection naturelle avec l'ensemble des parties $\{A'_1,A''_1\ldots, A'_s,A''_s\}$ de $Y_2$ associ\'ees aux translations $t_i$. Il nous suffit donc de montrer que la repr\'esentation du groupe  $T$ sur $Y_2$ d\'etermine toutes les partitions de $Y_2$ en $\{A'_i,A''_i\}$ pour en d\'eduire  l'unicit\'e du graphe $(\gL,Y)$.  Or, pour un indice $i$  et un point $u$ de $Y_2$ donn\'es dans $A'_i$, la condition $(5)$ de la proposition 7 nous montre que chaque point  $t_j(u)$ pour $j\neq i$ est dans $A''_i$ ainsi que les points $t_i\oo t_j(u)$ ce qui nous donne les $\ov{s}=2(s-1)$ points de $A''_i$. Ceux de $A'_i$ sont obtenus par exemple, par compl\'ementarit\'e. Ceci ach\`eve la preuve de l'unicit\'e.\\
L'existence se d\'emontre facilement par d\'enombrement puisqu'on peut repr\'esenter ''\`a la main'' les diff\'erents graphes  $Y_2$ :\\
$*$ Si $(t,s,\ov{s})=(1,2,2)$, $(\gL_2,Y_2)$ est un carr\'e.\\
$*$ Si $(t,s,\ov{s})=(1,3,4)$, $(\gL_2,Y_2)$ est un cube.\\
$*$ Si $(t,s,\ov{s})=(1,5,8)$, $(\gL_2,Y_2)$ est un hypercube en dimension $4$ auquel on rajoute ses diagonales principales.\\
Venons en \`a la d\'etermination des groupes d'automorphismes de ces graphes. On ne revient  pas sur le cas $s=1$.\\
Si $(t,s,\ov{s})=(1,2,2)$. On a $|Y|=1+2s+2\ov{s}=9$ et $|X|=|Y|+1=10$.\\
 Le groupe $\Aut(\gL_2,Y_2)$ est, par construction, isomorphe au groupe di\'edral $D_4$ des automorphismes d'un carr\'e et le groupe $\Aut(\gL,Y)$ \'etant transitif sur $Y$ il est d'ordre $9.|D_4|=72$. Enfin le groupe $\Aut(\gC,X)$ est \`a son tour transitif sur $X$, donc doublement transitif, et d'ordre $10.9.8=720$
 \\
Si  $(t,s,\ov{s})=(1,3,4)$. On a $|Y|=1+2s+2\ov{s}=15$ et $|X|=|Y|+1=16$.\\
 Le groupe $\Aut(\gL_2,Y_2)$ est, par construction, isomorphe au groupe $C$ des isom\'etries d'un cube d'ordre $|C|=48$,  et le groupe $\Aut(\gL,Y)$ \'etant transitif sur $Y$ il est d'ordre $15.|C|=720$. Enfin le groupe $\Aut(\gC,X)$ est \`a son tour transitif sur $X$, donc doublement transitif, et d'ordre $16.720=11520$.
 \\
Si  $(t,s,\ov{s})=(1,5,8)$.  On a $|Y|=1+2s+2\ov{s}=27$ et $|X|=|Y|+1=28$.\\
 Le groupe $\Aut(\gL_2,Y_2)$ est, par construction, isomorphe au
 groupe $H$ des isom\'etries d'un hypercube en dimension $4$
 compl\'et\'e par ses diagonales principales d'ordre
 $|H|=2^4.5!=1920$,  et le groupe $\Aut(\gL,Y)$ \'etant transitif sur
 $Y$ il est d'ordre $27.|H|=51840$. Enfin le groupe $\Aut(\gC,X)$ est
 \`a son tour transitif sur $X$, donc doublement transitif, et d'ordre
 $28.51840=1451520$, tout comme la repr\'esentation de degr\'e $28$ du
 groupe  $S_6(2)=PSp_6(2)$, ce qui permet de l'identifier. 
\hfill $\Box$ 
\begin{center}
  \includegraphics[scale=0.15]{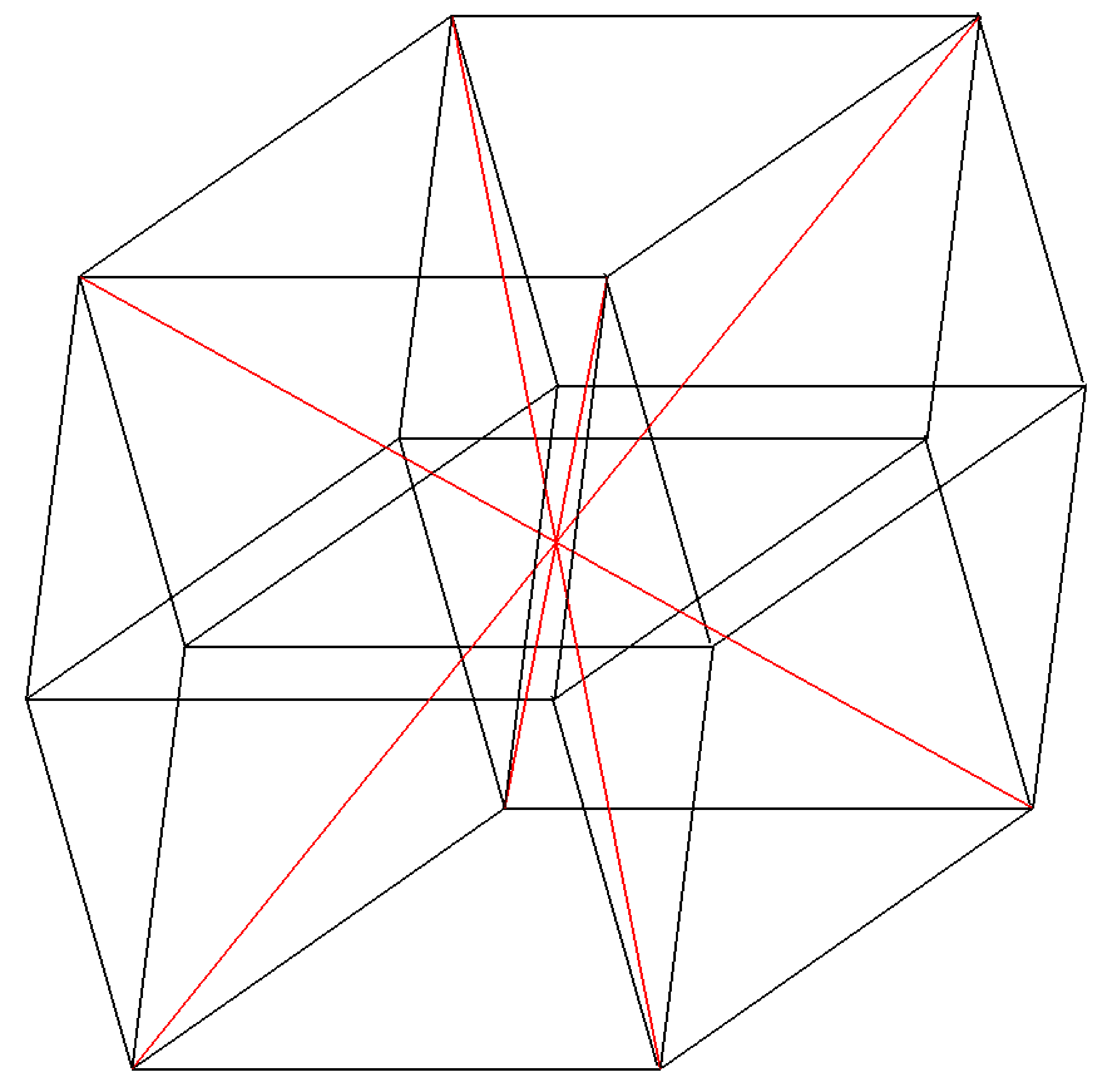}
\npb \centerline{L'hypercube et ses diagonales principales}
 \end{center}

 \subsection{Les graphes de Paley $P(q)$ o\`u $q=4t+5$}

Rappelons que si les param\`etres du graphe $(\gL,Y)$ sont $(t,s,\ov{s})$, ceux de son graphe compl\'ementaire $(\ov{\gL},Y)$ sont $(\ov{t},\ov{s},s)$  o\`u   $t+\ov{t}=s+\ov{s}-2$. On peut donc se contenter de rechercher l'un des deux graphes $(\gL,Y)$ ou $(\ov{\gL},Y)$.  Cette remarque nous pousse aussi \`a regarder en priorit\'e les graphes susceptibles d'\`etre isomorphes \`a leur compl\'ementaire, autrement dit pour lesquels  $(\ov{t},\ov{s},s)=(t,s,\ov{s})$. On a alors $s=\ov{s}=2s-t-1$, donc $s=\ov{s}=t+1$.  Or lorsque de plus $q=4t+5$ est une puissance d'un nombre premier $p$ il existe effectivement un graphe fortement r\'egulier avec ces param\`etres : le graphe de Paley $P(q)$ ([??]). On en donne ci-dessous une construction montrant que ce graphe est extensible.
  
On choisit un entier naturel $t$ tel que $q=4t+5$ soit une puissance d'un nombre premier $p$, et on pose $s=\ov{s}=t+1$. Soit  $F_q$ le corps fini \`a $q$ \'el\'ements,  $V=F_q\x F_q
$ et $X=\PPP(V)$  la droite projective associ\'ee qui est de cardinal
$q+1=4t+6$. 
Notons aussi $i$ une racine primitive quatri\`eme de $1$ dans $F_q$, $C$
l'ensemble des carr\'es non nuls et $\ov{C}$ l'ensemble des
non-carr\'es dans $F_q$. Enfin pour tout scalaire $\gl$, on pose 
$\gl+C=\{\gl +c\ |\ c\in C\}$. 
\begin{lemf}$\ $\\ 
$1$. $|C|=|\ov{C}|=(q-1)/2=2s$. \\
$2$. Si $\gl $ est un carr\'e non nul alors \ $|(\gl+C)\cap C|=s-1$ \  et  \   $|(\gl+C)\cap \ov{C}|= s$.\\
$3$. Si $\gm $ n'est pas un carr\'e alors \ $|(\gm+C)\cap C|=s$ \   et  \  $|(\gm+C)\cap \ov{C}|= s$.
\end{lemf}
\Dem \\
1. La premi\`ere affirmation est claire puisque l'application $x\to x^2$ de $F_q^*$ dans lui \mm est un morphisme d'image $C$ et de noyau $\{\pm 1\}$.\\
2. Pour montrer ces \'egalit\'es, on remarque tout d'abord que la multiplication par un carr\'e non nul $\gl$ induit une bijection de  $(1+C)\cap C$ sur $(\gl+C)\cap C$ et de  $(1+C)\cap \ov{C}$ sur $(\gl+C)\cap \ov{C}$, donc il nous suffit de les \'etablir lorsque $\gl=1$.\\
Recherchons le cardinal de l'ensemble $B$ des carr\'es non nuls $b^2$ aussi dans $1+C$. Dans le plan projectif, la conique projective d'\'equation $x^2+y^2=z^2$ comporte exactement $q+1$ points (voir $ [13]$ par exemple)  dont $6$ exactement ont une composante nulle, \`a savoir $\al(0,1,1)\ar$,$\al(0,1,-1)\ar$, $ \al(1,0,1)\ar$, $ \al(1,0,-1)\ar$, $  \al(1,i,0)\ar$ et $ \al(1,-i,0)\ar$. L'ensemble des $q-5$ points restant se projette surjectivement sur l'ensemble des carr\'es $b^2$ recherch\'es (par $\al(x,y,z)\ar \to z^2/x^2$), et chaque carr\'e $b^2$ poss\`ede exactement quatre ant\'ec\'edents par cette projection qui sont $\al(1,y,z)\ar$, $\al(1,-y,z)\ar$, $\al(1,y,-z)\ar$, $\al(1,-y,-z)\ar$. 
On en d\'eduit donc que le nombre des carr\'es $b^2$ qui sont aussi dans $1+C$ est $(q-5)/4=s-1$. 
La deuxi\`eme \'egalit\'e s'en d\'eduit facilement en remarquant que $\{C,
\ov{C}\}$ est une partition de $F_q^*$, que  $|C|=|\ov{C}|=|\gl+C|=2s$
\  et que \  $0\in (\gl+C)$ (car $(-1)$ est un carr\'e). On en tire  $|(\gl
+C)\cap (C\cup \ov{C})|=2s-1$ \  et \ $|(\gl+C)\cap
\ov{C}|=(2p-1)-(p-1)=p$.\\
3. Remarquant que la multiplication par un non-carr\'e
$\gm$  echange les ensembles $C$ et $\ov{C}$, on en d\'eduit que
$|(\gm+\ov{C})\cap C|=|(\gm^2+C)\cap\ov{C}|=p$,\  puis par translation
$|\ov{C}\cap (\gm+C)|=p$,\  enfin comme $\{C,
\ov{C}\}$ est une partition de $F_q^*$ \ et \ $(\gm+C)\subset F_q^*$,
\ il \  vient  $|C\cap (\gm+C)|=2p-p=p$. \hfill $\Box$\\

Donnons nous une base rang\'ee  $(u,v)$ de $V$ et l'application  \mbox{$\theta :\gl \to \al\gl.u +v \ar$} de $F_q$ dans l'ensemble  $X^u=\PPP(V)-\{\al u\ar\}$. Comme $\theta$ est clairement bijective on peut
associer \`a la base $(u,v)$ un graphe $\gC^u_{v}$ sur $X^u$ donn\'e  par ses ar\`etes : 
deux points $\al\ga.u+v\ar$ et $\al\gb.u+v\ar $ sont li\'es \si $ \ga-\gb \in C$, ce qu'on note aussi \\
\centerline{$\al\ga.u+v\ar\ \sim \ \al\gb.u+v\ar \  \ssi \ \ga-\gb \in C$}
(comme $(-1)$ est un carr\'e cette liaison  est sym\'etrique).\\
{\it Remarque.} On voit qu'en particulier le point $x=\al\ga.u+v\ar  $ est
li\'e \`a $\al v\ar  $ \si $\ga \in C$, mais aussi que deux points $x=\al \ga.u+v\ar  $ et
$y= \al \gb.u+v\ar  $  sont simultan\'ement li\'es ou non li\'es \`a $\al v\ar  $
\si $\ga\gb \in C$.\\ 
Le \ttt 3 nous montre que les graphes  $(\gC^u_{v},X^u)$, et nous d\'ecrit leur groupe d'automorphismes  ainsi que le groupe $G$ des permutations de $X$ qui permutent les graphes $(\gC^u_{v},X^u)$.

\begin{thmf}$\ $\\ 
On choisit une base $(u,v)$ de $V$.\\
$1$. Le diam\`etre du graphe $\gC^u_v $ est $2$. Pour deux points $x$
et $y$ dans $X^u$ on a\\
\hspace*{0.5cm} $\ast$ Si $d(x,y)=0$ alors $\gC^u_v(x,1)=\gC^u_v(y,1)$ est de cardinal $(q-1)/2=2s$.\\
\hspace*{0.5cm} $\ast$ Si $d(x,y)=1$ alors $\gC^u_v(x,1)\cap \gC^u_v(y,1)$ est de cardinal $t=s-1$.\\
\hspace*{0.5cm} $\ast$ Si $d(x,y)=2$ alors $\gC^u_v(x,1)\cap \gC^u_v(y,1)$ est de
cardinal $s$.\\
$2$. Le graphe $\gC^v_u $ est le localis\'e en \ $\al v\ar  $ \ de $\gC^u_v $.\\
$3.a.$  Pour tout automorphisme $\gvf \in GL_2(q)$ le graphe
$\gC^{\gvf(u)}_{\gvf(v)}$ sur $X^{\gvf(u)}$ est l'image par $\gvf$ du
graphe $\gC^u_v $ sur $X^u$ :
$\gvf(\gC^u_v )=\gC^{\gvf(u)}_{\gvf(v)}$.\\
$3.b.$  Supposons de plus que  $\gvf$ stabilise la droite $\al u\ar  $.\\
Si le d\'eterminant  $\det(\gvf) $  est un carr\'e dans le corps $F_q$, alors  $\gC^{\gvf(u)}_{\gvf(v)}=\gC^u_v $.\\
Sinon, le graphe 
$\gC^{\gvf(u)}_{\gvf(v)}$  est le graphe compl\'ementaire $\ov{\gC}^u_{v}$
de $\gC^u_v $ sur $X^u$.\\[0.1cm]
$4$. Le groupe $SL_2(V)$ poss\`ede deux orbites sur l'ensemble des
graphes $\gC^u_v $. L'orbite de l'un d'eux est l'ensemble
de ses localis\'es en tous les points de $X=\PPP(V)$.
\end{thmf}
\Dem \\
1. Posons $x=\al \ga.u+v\ar  $, $y=\al \gb.u+v\ar  $, $z=\al \gc.u+v\ar  $ et cherchons \`a quelle
condition $z$ est simultan\'ement li\'e \`a $x$ et
$y$. On a\\[0.1cm]
\vspace{0.2cm}\centerline{$x\sim z $ \ et \ $y\sim z $ \ $\ssi $ $\gc-\ga \in C$ \
  et \ $\gc-\gb \in C$ \ $\ssi $ $\gc \in (\ga +C)\cap(\gb +C)$ } 
Or par translation $|((\ga +C)\cap(\gb +C)|=|(\ga -\gb) +C)\cap
C|$. Utilisons le \nlb \mbox{lemme 4.}\\
$*$ Si $x=y$ alors \ $\ga=\gb$ \  et il vient $|\gC^u_v(x,1)|=|\gC^u_v(y,1)|=|C|=2s$.\\
$*$ Si $x$ et $y$ sont li\'es, alors \ $\ga-\gb$ \ est un
carr\'e non nul 
et d'apr\`es le lemme 4 \ $|((\ga -\gb) +C)\cap C|=s-1$. Il vient donc \vspace{0,1cm}$|\gC^u_v(x,1)|=|\gC^u_v(y,1)|=s-1$.\\
$*$ Enfin si $x$ et $y$ ne sont pas li\'es  alors \ $\ga-\gb\in \ov{C}$, et  \mbox{$|((\ga -\gb) +C)\cap C|=s>0$}. On en d\'eduit que $x$ et $y$ sont \`a distance $2$ et que  $|\gC^u_v(x,1)|=|\gC^u_v(y,1)|=s$.\\
2. Comparons les ar\`etes des graphes  $\gC^v_u $ et de  $\gC^u_v $.\\
Soient $x=\al \ga.u+v\ar  $ et  $y=\al \gb.u+v\ar  $ deux points de $\gC^u_v  $
distincts de $\al u\ar  $ et $\al v\ar  $ (donc $\ga\gb
\neq 0$). \\[0.1cm]
\hspace*{0.3cm} Dans le graphe $\gC^u_v  $ on a :\quad  $x=\al \ga.u+v\ar  \ \sim \ \al v\ar   $  $\ssi \ga \in
C$ , \\[0,1cm]
\hspace*{0.3cm} et dans le graphe $\gC^v_u  $ on a :\quad $x=\al u+(1/\ga)v\ar  \ \sim \ \al u\ar  $  $ \  \ssi \ 1/\ga \in C$.\\[0.1cm]
Donc les points \`a distance $1$ (resp. $2$) de $\al v\ar  $ dans  $\gC^u_v $ sont les
m\^emes que les
points \`a distance $1$ (resp. $2$) de $\al u\ar  $ dans $\gC^v_u $.\\
Par ailleurs $x$ et $y$ sont li\'es dans le graphe  $\gC^u_v $ \si
$\ga - \gb \in C$, tandis qu'ils sont li\'es dans le graphe
$\gC^v_u $ \si \vspace{0,1cm}$\displaystyle{\frac{1}{\ga}-\frac{1}{\gb}=\frac{\gb-\ga}{\ga\gb}\in C}$.
 Donc les
liaisons entre $x$ et $y$ sont conserv\'ees si $\ga\gb \in C$ et
invers\'ees dans le cas contraire. Mais d'apr\`es la remarque qui
pr\'ec\`ede le \ttt 3 cela signifie que la liaison entre deux points
$x$ et $y$ est invers\'ee ($x\sim y \lraa  x\not\sim y$) \  \si $x$ et
$y$ ne sont pas \`a la \mm distance de $\al v\ar  $ (ou $\al u\ar  $). 
Or ceci signifie pr\'ecis\'ement  que $\gC^v_u $ est le localis\'e de  $\gC^u_v $ (voir le paragraphe 2.2).\\
3.a. Chaque \elt $\gvf $ de $GL_2(V)$ agit naturellement sur l'ensemble des points
$\al w\ar  $ de $X=\PPP(V)$ par $\gvf(\al w\ar  )=\al \gvf(w)\ar  $, et en particulier
$\gvf(\al \ga.u+v\ar  )=\al \ga.\gvf(u)+\gvf(v)\ar  $, donc deux points
$x=\al \ga.u+v\ar  $ et $y=\al \gb.u+v\ar  $ de $X^u$ sont li\'es dans le graphe
$\gC^u_v $ \si leurs images \ $\gvf(x)$ \  et \ $\gvf(y)$ le sont dans le
graphe $\gC^{\gvf(u)}_{\gvf(v)}$.\\
3.b. Les graphes\  $\gC^u_v $ \
et \  $\gC^{\gvf(u)}_{\gvf(v)}$ ont \mm ensemble de sommets puisque  $X^u=X^{\gvf(u)}$.
L'image de la base  $(u,v)$ par $\gvf$ est de la forme $(a.u,
w=b.u+c.v)$ o\`u $ac\neq 0$, donc choisissant deux points $x=\al \ga a.u +w\ar  $ et $y=\al \gb
a.u +w\ar  $ de $X^u$, on voit qu'ils sont li\'es dans le graphe
$\gC^{\gvf(u)}_{\gvf(v)}$ \si $\ga -\gb \in C$, tandis qu'en les
r\'e-\'ecrivant $\displaystyle{x=\al (\frac{\ga a+b)}{c}.u+v\ar } $ et  $\displaystyle{y=\al (\frac{\gb a+b)}{c}.u+v\ar}  $ on
voit qu'ils sont li\'es dans le graphe $\gC^u_v $ \si $(\ga-\gb)a/c
\in C$. Ainsi les liaisons des graphes $\gC^u_v $ et
$\gC^{\gvf(u)}_{\gvf(v)}$ sont toutes identiques lorsque $a/c \in C$ et
toutes diff\'erentes lorsque $a/c \in \ov{C}.$ Comme \  $\det
\gvf=(a/c).c^2$,\  le r\'esultat en d\'ecoule.\\
4. D'apr\`es $3a$ le groupe $PGL_2(V)$ op\`ere transitivement sur
l'ensemble $\gC$ des graphes $\gC^u_v $. Donc, comme le groupe $PSL_2(V)$ est
d'indice 2 dans $PGL_2(V)$,  soit il est lui \mm transitif sur $\gC$,
soit il y poss\`ede deux orbites. Mais prenant l'un des graphes
$\gC^u_v $, on voit qu'il n'existe aucun \elt $\gvf$ de $SL_2(V)$ tel que
$\gvf(\gC^u_v )=\gC^{\gvf(u)}_{\gvf(v)}=\ov{\gC}^u_{v}$, car cela
impliquerait tout d'abord que $\gvf$ stabilise la droite vectorielle
$\al u\ar  $ puis, gr\^ace \`a $3a.$, que
$\gC^{\gvf(u)}_{\gvf(v)}=\gC^u_v =\ov{\gC}^u_{v}$.
Donc $SL_2(V)$ poss\`ede deux orbites sur $\gC$, et pour terminer il nous
faut montrer que tout localis\'e d'un graphe $\gC^u_v $ est dans
l'orbite de  $\gC^u_v $ sous l'action de $SL_2(V)$. Mais ceci vient
du fait (vu en $2$) 
que le localis\'e de $\gC^u_v $ en $\al v\ar  $ est  $\gC^v_u $ et que le
morphisme $\gvf$ qui envoie la base $(u,v)$ sur la base $(v,u)$ admet
pour d\'eterminant $(-1)$ qui est un carr\'e dans $F_q$. \hfill $\Box$

\section{Repr\'esentations r\'eduites des graphes} 
Dans cette partie on associe aux graphes pr\'ec\'edemment calcul\'es leurs repr\'esen\-tations r\'eduites. Le principe repose sur l'article ([5]) et plus particuli\`erement sur ses propositions 1, 2 et le \ttt 1. On y prouve que toute repr\'esentation r\'eduite $u:X\to (E,\ga)$ d'un graphe $(\gC,X)$  est caract\'eris\'ee par ses param\`etres $(\go,c)$ et donc par sa matrice $S(u)=S(\go,c)$. 
Pour $X=\{1,\ldots,n\}$, soit $x_1,\ldots,x_n$ une base orthonom\'ee d'un espace euclidien $E$ dont le produit scalaire est not\'e $\gvf=(\ | \ )$ et $x : X\to E$ la repr\'esentation triviale de $\gC$ donn\'ee par $x : i\to x_i$ pour tout $i$ dans $X$. La forme bilin\'eaire sur $E$ dont la matrice dans la base $x_1,\ldots,x_n$  est $\EE=S(1,1)$  est not\'ee $\gb$ ainsi que l'endomorphisme sym\'etrique qui lui est associ\'e par  \\[0,1cm]
\vspace{0,1cm}\centerline{$\forall y,z \in E,\quad \gb(x,y)=(\gb(x)|y)=(x|\gb(y)).$}
\`A chaque nombre r\'eel $\gl$ on associe la forme bilin\'eaire $\gb-\gl.Id$ de matrice $S(1-\gl,1)$ et de rang $n-m(\gl)$ o\`u $m(\gl)$  est la multiplicit\'e de $\gl$ comme valeur propre de $S(1,1)$ (\'eventuellement nulle). Soit enfin $p_\gl$ le projecteur orthogonal de $E$ sur l'image $I_\gl$ de $\gb-\gl.Id$ et $E_\gl=\ker p_\gl$ son noyau. Le rang de $S(1-\gl,1)$ \'etant \'egal \`a la dimension de $I_\gl$, la repr\'esentation $u_\gl:=p_\gl\oo x$ de $\gC$ sur $I_\gl$ est r\'eduite, de matrice $S(1-\gl,1)$, ce qui nous donne une construction  d'une repr\'esentation r\'eduite de matrice donn\'ee. \\
Pour calculer explicitement les repr\'esentations r\'eduites des graphes \'etudi\'es dans la section 3 il nous faut bien s\^ur calculer les valeurs propres de la matrice $S_\gC(1,1)$ de chacun d'eux. Nous allons tout d'abord \'etablir un r\'esultat int\'eressant en soi. 

 \subsection{Nombre de valeurs propres de la matrice $\EE=S(1,1)$} 

\begin{propf}$\ $\\ 
Soit $(\gC, X)$ un graphe de matrice $\EE=S(1,1)$.\\
Si le groupe $G(\gC)$ op\`ere deux fois transitivement sur l'ensemble $X$, alors la matrice $S(1,1)$ poss\`ede au plus deux valeurs propres distinctes.
\end{propf}
\Dem. 
On garde les notations ci-dessus. Soit $\gl_1,\ldots,\gl_s$ les valeurs propres de $S(1,1)$. L'espace $E$ est une somme orthognale des sous-espaces propres $E_t:=E_{\gl_t}$ de la forme $\gb$ de matrice $S(1,1)$ $(1\leq t \leq s)$ : \\[0,1cm]
\vspace{0,1cm}\centerline{ $E =E_1\oplus\cdots\oplus E_s$}
Par d\'efinition du groupe $G=G(\gC)$, si $\gs$ est dans $G$, il existe une suite de coefficients $\nu_1,\ldots,\nu_n$ dans $\{-1,+1\}$ tels que \\[0,1cm]
$(3)$ \vspace{0,2cm}\centerline{\hspace{3cm}$\forall (i, j)\in X\x X ,
  \qquad \gve_{\gs(i),\gs(j)}=\nu_i\nu_j.\gve_{i,j}$ \hspace{3cm},}
ce qui nous permet d'associer \`a $\gs$ une $\gb$-isom\'etrie $f_\gs$ de $E$ en posant  \\[0,1cm]
 \vspace{0,2cm}\centerline{$\forall i\in X,\quad f_\gs(x_i)=\nu_i.x_{\gs(i)}$.}
 Mais $f_\gs$ est aussi une $\gvf$-isom\'etrie de $E$ puisque l'image par $f_\gs$ de $(x_1,\ldots,x_n)$ est une base orthonorm\'ee de $(E,\gvf)$.  On en d\'eduit donc que pour deux points $y$ et $z$ quelconques de $E$ on a \\[0,1cm]
 \vspace{0,2cm}\centerline{$ (f_\gs\oo\beta(y) | z)=(\beta(y) | f^{-1}_\gs (z))=\beta(y, f^{-1}_\gs (z))=\beta(f_\gs (y),  z)=(\beta \oo f_\gs( y),  z)$.}
 Par cons\'equent $f_\gs\oo\beta=\beta \oo f_\gs$. On en d\'eduit que les sous-espaces propres $E_t$ de $\beta$ sont stables par l'isom\'etrie $f_\gs$, et si, pour tout vecteur $y$ et tout indice $t$ entre $1$ et $s$, on note $q_t(y):=y^t$ la projection orthogonale de $y$ sur $E_t$, on a \\[0,1cm]
 \vspace{0,2cm}\centerline{ $\forall y\in E,\quad  f_\gs(q_t(y))=q_t(f_\gs(y))$.}
 Montrons que l'application $u_t:=q_t\oo x : X\to E_t$ est une repr\'esentation de $X$.\\
Prenons quatre points $i,j,k,l$ de $X$ tels que $i\neq j$ et $k\neq l$. Le groupe $G$ op\'erant doublement transitivement sur $X$, il existe une permutation $\gs$ de $X$ telle que $\gs(i)=k$ et $\gs(j)=l$. D\'efinissons des constantes $c_t$ et $\go_t$ par les egalit\'es \\[0,1cm]
 \vspace{0,2cm}\centerline{ $(x_i^t|x_j^t)=\gve_{i,j}c_t$\quad et \quad $(x_i^t|x_i^t)=\go_t$. }
 Il vient en remarquant que \ $\gve_{k,l}=\gve_{\gs(i),\gs(j)}=\nu_i\nu_j.\gve_{i,j}$, \\[0,2cm]
 \vspace{0,2cm}\centerline{$(x_k^t|x_l^t)=(\nu_i.f_\gs(x_i^t)|\nu_j.f_\gs(x_j^t))= \nu_i\nu_j(x_i^t|x_j^t)= \nu_i\nu_j\gve_{i,j}c_t=\gve_{k,l}c_t$.}
 Par cons\'equent la constante $c_t$ est ind\'ependante du couple d'indices distincts $(i,j)$. Un raisonnement similaire montre que la constante $\go_t$ ne d\'epend pas de l'indice $i$ et ceci montre que l'application $u_t:=q_t\oo x : X\to E_t$ est une repr\'esentation de $X$ sur $E_t$, de matrice $S(\go_t,c_t)$. On a rappel\'e plus haut que les repr\'esentations r\'eduites de $\gC$ sont de degr\'e $r_t=n-m(\gl_t)$ o\`u $\gl_1,\ldots,\gl_s$ d\'esignent les valeurs propres de la matrice $S(1,1)$. Par ailleurs, d'apr\`es la proposition 2 de ([5]) toute repr\'esentation d'un graphe est somme d'une repr\'esentation nulle et d'une repr\'esentation r\'eduite. La dimension  $m(\gl_t)$ de chaque sous-espace propre $E_t$ de $S(1,1)$ doit donc \^etre  sup\'erieure ou \'egale \`a l'un des nombres $r_j=n-m(\gl_j)$.
 Comme $n=m(\gl_1)+\cdots+m(\gl_s)$ on ne peut qu'avoir pour tout $t$, $m(\gl_t)\geq n-m(\gl_t) $. En additionnant ces in\'egalit\'es il vient $n\geq s.n-n$ donc $1\geq s-1$ soit encore $2\geq s$ comme attendu. \hfill $\Box$

{\it Remarque}. Lorsque $\gl\neq 1$ les repr\'esentations de matrices $S(1-\gl,1)$ et $S(1,1/(1-\gl))$ sont homoth\'etiquement semblables et le coefficient $c=1/(1-\gl)$ est g\'eom\'etriquement int\'eressant puisqu'il repr\'esente le cosinus de l'angle commun de deux droites dans la gerbe de la repr\'esentation.

 \subsection{Application aux graphes} 
Dans la section 3 nous avons obtenu diff\'erents graphes extensibles $(\gC^x, X^x)$, de param\`etres $(t,s,\ov{s})$. Chacun d'eux est le localis\'e  d'un graphe $(\gC,X)$ dont les repr\'esentations r\'eduites  rang $r<n=|X|$ sont associ\'ees, d'apr\`es ce qui pr\'ec\`ede, soit \`a la matrice $S(0,1)$ lorsqu'elle est  singuli\`ere, soit aux racines du polyn\^ome $\chi(x)=\det(S(1,x))$. En utilisant Maple, on constate que parmi ces exemples, aucune des matrices $S(0,1)$ n'est singuli\`ere. On calcule alors dans chaque cas 
le polyn\^ome $\chi(x)=\det(S(1,x))$. 

$*$ Graphe $(\gC^x, X^x)$ :  $(t,s,\ov{s})=(0,1,1)$. On a  $|X|=2+2s+2\ov{s}=6$,\\[0,2cm]
\vspace{0,2cm}\centerline{$\chi(c)=-(5c^2-1)^3$ \quad et \quad  $\displaystyle{\cos\theta=\pm \frac{1}{\sqrt{5}}}$ }
Le graphe $(\gC,X)$ admet donc une repr\'esentation $u$ dans un espace euclidien $E$ de dimension $3$. La gerbe de droites $\GG(u)=\{\al u_1\ar,\ldots,\al u_6\ar\}$  compte six droites formant deux \`a deux un angle $\theta$ dont le cosinus vaut \vspace{0,1cm}$\pm 1/\sqrt{5}$. On reconnait ici les six axes passant par les centres des faces oppos\'ees d'un dod\'eca\`edre.\\
Le groupe d'automorphisme $G(u)$ est semblable au groupe altern\'e $A_5$ agissant doublement transitivement sur un ensemble de cardinal 6.

$*$ Graphe $(\gC^x, {X^x})$ :  $(t,s,\ov{s})=(1,1,0)$. On a   $|X|=2+2s+2\ov{s}=4$\\[0,2cm]
\vspace{0,2cm}\centerline{ $\chi(c)=-(3c-1)(c+1)^3$ \quad et \quad  $\displaystyle{\cos\theta=\frac{1}{3}}$ \  ou \  $\displaystyle{\cos\theta=-1}$}
Le graphe $(^x\gC, {^xX})$ est un ''triangle point\'e '', \cad un triangle auquel on a rajout\'e un point isol\'e. C'est aussi le localis\'e en $x$ d'un t\'etra\`edre. Il admet deux repr\'esentations r\'eduites correspondant \`a ces deux  pr\'esentations. La premi\`ere est donn\'ee par la gerbe des quatre droites   joignant le centre d'un t\'etra\`edre r\'egulier \`a ses sommets ($\textstyle{\cos\theta=\frac{1}{\sqrt{3}}}$). La seconde ($\textstyle{\cos\theta=-1}$) envoie les trois sommets du triangle de  $(^x\gC, {^xX})$  sur le point $1$ de l'axe r\'eel et le point isol\'e $x$ sur $-1$. 
 
$*$ Graphe $(\gC^x, X^x)$ :  $(t,s,\ov{s})=(1,3,4)$. On a   $|X|=2+2s+2\ov{s}=16$ \\[0,2cm]
\vspace{0,2cm}\centerline{ $\chi(c)=(5c+1)^6(3c-1)^{10}$ \quad et \quad  $\displaystyle{\cos\theta=-\frac{1}{5}}$ \  ou \  $\displaystyle{\cos\theta=\frac{1}{3}}$ }
Le graphe $(\gC,X)$ admet deux repr\'esentations r\'eduites de degr\'e $<16$ sur une gerbe isom\'etrique de 16 droites. 
L'une dans un espace euclidien de dimension 6, les droites faisant deux  \`a deux un angle $\theta$ tel que $\cos\theta =1/3$, l'autre dans un dans un espace euclidien de dimension 10, les droites faisant deux  \`a deux un angle $\theta$ tel que $\cos\theta =-1/5$.\\
Le groupe $G(u)$ des isom\'etries stabilisant ces gerbes de droites est semblable au groupe $G(\gC)$ agissant sur $X$. Il est d'ordre  $11520$  et op\`ere deux fois transitivement sur la gerbe $\GG(u)$.

$*$ Graphe $(\gC^x, X^x)$ :  $(t,s,\ov{s})=(1,5,8)$. On a   $|X|=2+2s+2\ov{s}=28$\\[0,2cm]
\vspace{0,2cm}\centerline{  $\chi(c)=-(9c+1)^7(3c-1)^{21}$ \quad et \quad $\displaystyle{\cos\theta=-\frac{1}{9}}$ \  ou \ $\displaystyle{\cos\theta=\frac{1}{3}}$}
Le graphe $(\gC,X)$ admet deux repr\'esentations r\'eduites de degr\'e $<28$ sur une gerbe isom\'etrique de 28 droites. 
L'une dans un espace euclidien de dimension 7, les droites faisant deux  \`a deux un angle $\theta$ tel que $\cos\theta =1/3$, l'autre dans un dans un espace euclidien de dimension 21, les droites faisant deux  \`a deux un angle $\theta$ tel que $\cos\theta =-1/9$.\\
Le groupe $G(u)$ des isom\'etries stabilisant ces gerbes de droites est semblable au groupe $G(\gC)$ agissant sur $X$. Il est d'ordre $1451520=2.725760$  et op\`ere deux fois transitivement sur la gerbe $\GG(u)$.

$*$ Graphe $(\gC^x, X^x)$ : $(t,s,\ov{s})=(t,t+1,t+1)$.   $|X|=2+2s+2\ov{s}=q+1$. \\
Le graphe $(\gC^x,X^x)$ est un graphe de Paley $P(q)$,  pour $q=4t+5$.\\[0,2cm]
\vspace{0,2cm}\centerline{ $\displaystyle{\chi(c)=\gve(c\sqrt{q}-1)^\frac{q+1}{2}(c\sqrt{q}+1)^\frac{q+1}{2}}$, o\`u $\gve=\pm 1$ \  et \  $\displaystyle{\cos\theta=\pm \frac{1}{\sqrt{q}}}$ }
Le graphe $(\gC,X)$ admet deux repr\'esentations r\'eduites de degr\'e $<q+1$ sur une gerbe isom\'etrique de q+1 droites plac\'ees dan un espace de dimension $q+1)/2$. Elles forment deux \`a deux un angle $\theta$ dont le cosinus vaut  $\pm 1/\sqrt{q}$.\\
Le groupe $G(u)$ des isom\'etries stabilisant ces gerbes de droites
est semblable au groupe $PSL_2(q)$  agissant sur l'ensemble $X$ des
points de la droite projective $\PPP_1(F_q)$ qui est d'ordre $(q+1)q(q-1)/2$.  

\section{Prolongements} 
Ce travail n'est bien s\^ur pas ferm\'e.
Tout d'abord il reste \`a compl\'eter la liste des graphes $(\gL,Y)$ extensibles.  Comme nous l'avons mentionn\'e plus haut ce sont tous des graphes fortement r\'eguliers dont les param\`etres, habituellement  not\'es $(v,k,\gl,\mu)$  sont reli\'es aux notations de cet article par\\[0,1cm]
 \vspace{0,1cm}\centerline{ $v=n$, $k=2s$, $\gl=t$, $\mu=s$.}
Mais inversement, si les param\`etres $(v,k,\gl,\mu)$ d'un graphe fortement r\'egulier $(\gL,Y)$ satisfont \`a ces relations pour des entiers $ n,s,t$ convenablement choisis, est on s\^ur que le graphe $(\gL,Y)$ est extensible ?\\
En particulier il existe des graphes $(\gL,Y)$ de param\`etres \mbox{$(t,t+1,t+1)$} (i.e. $v=4t+5$, $k=2t+2$, $\gl=t$, $\mu=t+1$), dit "conference graphs", qui ne sont pas des graphes de Paley. Lesquels sont extensibles ? Quels groupes leurs sont associ\'es ?

\section{Bibliographie}

\label{thebib}

{\it Articles}\\

[1] \quad Dominique de Caen, {\it Large equiangular sets of lines in Euclidean spaces.} The electronic journal of comninatorics. November 9, 2000.

[2] \quad L. Nguyen Van The, {\it On a problem of Specker about Euclidean
representations of finite graphs.}, (2008)\\
http://arxiv.org/abs/0810.2359

[3] \quad Aidan Roy , {\it Minimal Euclidean representations of graphs}. (2008)\\
 http://arxiv.org/abs/0812.3707

[4] \quad J.J. Seidel, {\it Discrete Non-Euclidean Geometry,} pp. 843-920 in Handbook of Incidence Geometry (F. Buckenhout, ed.), Elsevier 1995.

[5] \quad L. Vienne, {\it Repr\'esentations lin\'eaires des graphes finis.} (2009)\\ 
	http://arxiv.org/abs/0902.1874.

[6] \quad Blokhuis, A. {\it On Subsets of $GF(q^2)$ with Square Differences.} Indag. Math. 46, 369-372, (1984). 	

[7] \quad  Brouwer, A. E.; Cohen, A. M.; and Neumaier, A. {\it Conference Matrices and Paley Graphs.} In Distance Regular Graphs. New York: Springer-Verlag.

{\it Livres}

$[8]$ \quad Godsil, Chris, and Royle, {\it  Algebraic Graph Theory.} \mbox{New York: Springer. (2001)} 

$[9]$\quad  Dembowski, P. {\it   Finite  geometries.} Springer-Verlag (1968)

$[10]$\quad  Norman Biggs. {\it Finite groups of automorphisms.} \mbox{London Mathematical Society. (1970)}

$[11]$\quad  Norman Biggs. {\it Algebraic Graph Theory.} \mbox{Cambridge University Press,  (1993)}

$[12]$\quad P. J. Cameron and J.H.Van Lint. {\it Designs, graphs,codes and their links, vol 22 of London Mathematical Society Student Texts,} Cambridge University Press, Cambridge, (1991).

$[13]$ \quad L. Vienne. {\it Pr\'esentation alg\'ebrique de la g\'eom\'etrie classique}. \mbox{Vuibert. (1996).}

\end{document}